\documentclass{amsart}
\usepackage{amssymb,amsfonts}
\usepackage{enumerate}
\usepackage{mathrsfs}
\usepackage{graphicx}
\usepackage{mathtools}
\usepackage{leftidx}
\usepackage{tikz-cd}
\usepackage{bm}
\usepackage{url}
\newtheorem{thm}{Theorem}[section]
\newtheorem{cor}[thm]{Corollary}
\newtheorem{prop}[thm]{Proposition}
\newtheorem{lem}[thm]{Lemma}

\theoremstyle{definition}
\newtheorem{defn}[thm]{Definition}

\newtheorem{exmp}[thm]{Example}

\theoremstyle{remark}
\newtheorem{rem}[thm]{Remark}

\makeatletter
\let\c@equation\c@thm
\makeatother
\numberwithin{equation}{section}

\def\bthm{\begin{thm}}
\def\ethm{\end{thm}}
\def\blm{\begin{lem}}
\def\elm{\end{lem}}
\def\bdf{\begin{defn}}
\def\edf{\end{defn}}
\def\bpf{\begin{proof}}
\def\epf{\end{proof}}
\def\bpp{\begin{prop}}
\def\epp{\end{prop}}
\def\bcor{\begin{cor}}
\def\ecor{\end{cor}}
\def\brm{\begin{rem}}
\def\erm{\end{rem}}
\def\beg{\begin{exmp}}
\def\eeg{\end{exmp}}
\def\bL{\mathbb{L}}

\def\bX{\mathbb{X}}

\def\bZ{\mathbb{Z}}
\def\cC{\mathcal{C}}

\def\cG{\mathcal{G}}

\def\cO{\mathcal{O}}
\def\cP{\mathcal{P}}
\def\cQ{\mathcal{Q}}
\def\cR{\mathcal{R}}

\def\cY{\mathcal{Y}}

\def\frX{\mathfrak{X}}

\def\frg{\mathfrak{g}}
\def\frh{\mathfrak{h}}

\def\frm{\mathfrak{m}}

\newcommand{\raq}{\,\rightarrow \,}

\newcommand{\rontoq}{\,\twoheadrightarrow\,}

\newcommand{\xraq}[2][]{\, \xrightarrow[#1]{#2} \,}
\newcommand{\xlaq}[2][]{\, \xleftarrow[#1]{#2} \,}

\newcommand{\ra}{\rightarrow}
\newcommand{\la}{\leftarrow}
\newcommand{\rinto}{\hookrightarrow}

\newcommand{\ronto}{\twoheadrightarrow}

\newcommand{\xra}[2][]{\xrightarrow[#1]{#2}}
\newcommand{\xla}[2][]{\xleftarrow[#1]{#2}}

\newcommand{\ie}{{\it i.e.}}
\newcommand{\eg}{{\it e.g.}}
\newcommand{\cf}{{\it cf.}}

\newcommand{\Mod}{{\rm Mod}}

\newcommand{\Alg}{{\rm Alg}}

\newcommand{\Sym}{{\rm Sym}}

\newcommand{\Ch}{{\rm Ch}}
\newcommand{\DGA}{{\rm DGA}}
\newcommand{\DGLA}{{\rm DGLA}}
\newcommand{\CDGA}{{\rm CDGA}}
\newcommand{\dgMod}{{\rm dgMod}}

\newcommand{\DGAp}{{\rm DGA}^{\geq 0}}

\newcommand{\Tot}{{\rm Tot}}

\newcommand{\cone}{{\rm cone}}

\newcommand{\CE}{{\rm CE}}

\newcommand{\Ob}{{\rm Ob}}
\newcommand{\op}{{\rm op}}

\newcommand{\id}{{\rm id}}

\newcommand{\Hom}{{\rm Hom}}
\newcommand{\Homcom}{\underline{{\rm Hom}}}

\newcommand{\cHom}{\mathscr{H}\text{\kern -3pt {\calligra\large om}}\,}

\newcommand{\Map}{{\rm Map}}
\newcommand{\holim}{{\rm holim}}

\newcommand{\hocolim}{{\rm hocolim}}

\newcommand{\Spec}{{\rm Spec}}

\newcommand{\dAff}{{\rm dAff}}

\newcommand{\Der}{{\rm Der}}

\newcommand{\gr}{{\rm gr}}

\newcommand{\OCar}{\Omega_{{\rm Car}}}

\newcommand{\DR}{{\rm DR}}

\newcommand{\DRcl}{{{\rm DR}}^{{\rm cl}}}

\newcommand{\DRCar}{{{\rm DR}}_{{\rm Car}}}

\newcommand{\DRCarcl}{{{\rm DR}}_{{\rm Car}}^{{\rm cl}}}

\newcommand{\XCar}{\frX_{{\rm Car}}}

\newcommand{\PCar}{\Pol_{{\rm Car}}}

\newcommand{\hPCar}{\widehat{\Pol}_{{\rm Car}}}

\newcommand{\CDGAn}{{\rm CDGA}^{\leq 0}}

\newcommand{\CDGAp}{{\rm CDGA}^{\geq 0}}

\newcommand{\setdel}{{\rm Set}_{\Delta}}

\newcommand{\AW}{{\rm AW}}

\newcommand{\CAlg}{{\rm CAlg}}

\newcommand{\Pol}{{\rm Pol}}

\newcommand{\CoSym}{{\rm CoSym}}

\newcommand{\hG}{{\rm h}G}

\newcommand{\hc}{{\rm hc}}

\newcommand{\dgModcom}{\underline{\rm dgMod}}

\newcommand{\El}{{\rm El}}

\newcommand{\Ran}{{\rm Ran}}

\newcommand{\hca}{{\rm hca}}

\newcommand{\ca}{{\rm ca}}

\newcommand{\Ad}{{\rm Ad}}

\newcommand{\CG}{{\rm CG}}

\newcommand{\vc}{{\rm vc}}

\DeclareMathAlphabet{\mathpzc}{OT1}{pzc}{m}{it}

\title{Shifted symplectic and Poisson structures on global quotients}
\author{Wai-Kit Yeung}
\address{Kavli IPMU, The University of Tokyo}
\email{wai-kit.yeung@ipmu.jp}
\begin{document}

\begin{abstract}
For a derived stack obtained as a quotient of a derived affine scheme by a reductive group, we show that shifted symplectic structures can be characterized by the Cartan-de Rham complex. For non-reductive groups, we also show the analogous statement for Getzler's extension of the Cartan-de Rham complex. Dually, we construct a Cartan model for polyvector fields on global quotients by reductive groups, and show that shifted Poisson structures can be characterized by it.
\end{abstract}

\maketitle


\tableofcontents

\section{Introduction}  \label{sec_intro}

In the papers \cite{PTVV, CPTVV} (see also \cite{Pri17}), the notions of shifted symplectic and Poisson structures were introduced. It has since then attracted the interest of many mathematicians, as such structures arise naturally in many mathematical contexts, especially when moduli spaces are concerned.

While the papers \cite{PTVV, CPTVV, Pri17} have set up the foundations for the theory, most of the available examples are given in terms of some general principles, mostly arising from the AKSZ construction. It seems that the only example in {\it loc. cit.} with an explicit%
\footnote{Of course, mathematicians disagree vastly on what it means by ``explicit''. We try to demonstrate what we mean in this introduction.} 
description is $BG = [*/G]$ for a reductive group $G$. 
In this paper, we work out an explicit description of shifted symplectic and Poisson structures for global quotients, \ie, derived stacks of the form $[Y/G]$, where $Y$ is an affine derived scheme.
These examples arise naturally when one considers, \eg, moduli spaces of representations of associative algebras, moduli spaces of local systems, {\it etc.}, all of which are classical subjects, with vast literature.
Our description will involve classical constructions such as the Cartan model for equivariant cohomology, and is therefore more congenial to these literature. 
We also expect our description to be useful in questions about quantization, for which shifted Poisson structures are expected to be a natural setting. Indeed, quantizations in representation theory are often obtained by writing explicit formulas (often guided by considerations of Hall algebras, convolution algebra, {\it etc.}). Our description of shifted symplectic and Poisson structures should be more congenial to these considerations, and could be useful, even if only as a vanguard, for the discovery of unsuspected linkages.
Another advantage is that our description of equivariant differential forms and polyvectors make heuristic sense for certain infinite dimensional (gauge) group $\cG$ acting on infinite dimensional derived spaces $\cY$, in which case we expect our description to serve as a useful heuristic guide.

Thus, consider the global quotient $X = [Y/G]$ for $Y = \Spec \, A$ and $G$ reductive. Let us try to guess what is a correct notion of shifted symplectic and Poisson structures. Just like in the classical case, shifted symplectic structures should depend on a certain de Rham complex
\begin{equation}  \label{DR_complex_intro}
\DR^0(X) \xraq{d} \DR^1(X) \xraq{d} \DR^2(X) \xraq{d} \ldots
\end{equation}
Here $\DR^p(X)$ should be the (derived) global sections of $\Lambda^p \bL_{X/k}$ where $\bL_{X/k}$ is the cotangent complex of $X$. As such, each $\DR^p(X)$ is not just a vector space, but a cochain complex $(\DR^p(X),\partial)$. Accordingly, \eqref{DR_complex_intro} is a bicomplex.

Shifted symplectic structures are defined in terms of closed $2$-forms. Thus, we should consider $\omega_2 \in \DR^2(X)$ that is closed under $d$. But in the context of derived stacks, since $\DR^3(X)$ is a cochain complex, it is not natural to require that $d(\omega_2) = 0$. Instead, one should require $d(\omega_2) = \partial(\omega_3)$ for some $\omega_3$. One should also want some control over $\omega_3$, say by requiring that $d(\omega_3) = \partial(\omega_4)$, and so on. In other words, a natural definition%
\footnote{There are also more sopihsticated reasons for making this definition in terms of the $S^1$-action on the derived loop stack. See, \eg, \cite{PTVV}.}
for a closed form would be a cocycle in the direct product total complex of $\DR^{\geq 2}(X)$.

In the case $X = [Y/G]$ of global quotient, recall that every $G$-equivariant DG module over $A$ gives rise to a quasi-coherent sheaf on $X$ (see Appendix \ref{app_cart} for a more detailed description). In particular, the cotangent complex is represented%
\footnote{In this introduction, we assume that $A$ is cofibrant, or more generally almost cofibrant in the sense of Definition \ref{almost_cof_def} below. In general, one should first $G$-equivariantly resolve the given commutative algebra $A$ before proceeding with our present considerations.} by the $G$-equivariant DG module
\begin{equation}  \label{OCar_intro}
\OCar^1(Y/G) \, = \, \cone[\, \Omega^1(A) \xraq{\alpha} \frg^* \otimes A \,][-1]
\end{equation}
where $\alpha$ is the $A$-linear map dual to the infinitesimal action $\frg \ra \Der(A)$. Also, since $G$ is reductive, taking global sections corresponds to taking $G$-invariants of the $G$-equivariant DG module ({\it cf.} Appendix \ref{app_cart}). This gives a description of $\DR^n(Y/G)$ as a cochain complex:
\begin{equation}  \label{DR_n_intro}
\DR^n(Y/G)[-n] \,= \, (\Sym^n_A( \OCar^1(A)[-1]))^G \, = \, \bigoplus_{p+q = n} \left( \, \Omega^p(A)[-p] \otimes \Sym^q(\frg^*[-2]) \, \right)^G
\end{equation}
with differentials induced from the cone \eqref{OCar_intro}.

It is a classical fact that, if we take $d' = d \otimes \id$ on the right hand side of \eqref{DR_n_intro}, then the resulting map $d' : \DR^n(Y/G) \ra \DR^{n+1}(Y/G)$ becomes a bicomplex (see Theorem \ref{Omega_CE_dR_inv}). Indeed, this is precisely the bicomplex appearing in the Cartan model for equivariant cohomology. We call it the Cartan-de Rham bicomplex, and will write it as $\DRCar^{\bullet}(Y/G)$ from now on.

This allows us to define a closed $2$-form (also called a pre-symplectic structure) on $X$. Precisely, an $m$-shfited pre-symplectic structure on $X$ is a cocycle of degree $m$ in a direct product total complex:
\begin{equation*}
\tilde{\omega} = (\omega_2, \omega_3,\ldots) \, \in \,  Z^m \biggl( \, \prod_{p \geq 2} \, \DRCar^p(Y/G)[2-p] \, , \,  \partial + d'  \, \biggr)
\end{equation*}

In particular, every $m$-shifted closed $2$-form has an underlying $m$-shifted $2$-form $\omega_2 \in Z^m(\DRCar^2(Y/G))$. Recall that $\DRCar^2(Y/G) := (\Lambda^2_A \,\OCar^1(A))^G$ (we write $\Lambda^p_A(M) := \Sym_A^p(M[-1])[p]$). Thus, in particular, $\omega_2$ induces a $G$-equivariant%
\footnote{There is some minor technical problem for giving a $G$-equivariant structure to the dual of a $G$-equivariant module, but it won't not concern us here. See Section \ref{Cartan_sec} for more details.} map of DG modules over $A$
\begin{equation*}
\omega_2^{\sharp} \, : \, \OCar^1(A)^{\vee}[-m] \raq \OCar^1(A)
\end{equation*}
The $m$-shfited pre-symplectic structure $\tilde{\omega}$ is said to be an $m$-shfited symplectic structure in Cartan model ({\it cf.} Definition \ref{shifted_sympl_def}) if this map is a quasi-isomorphism.

Dually, shifted Poisson structures depend on polyvector fields. Morally speaking, one would like to define the complex of $p$-polyvectors as the derived global sections of  $\Lambda^p((\bL_{X/k}^1)^{\vee})$.
However, in general, unless a module is finitely generated and projective, there might be pathologies in taking $\Lambda^p$ of its dual. Instead, one should first take $\Lambda^p$ and then take dual. 
Another modification is that we should take a certain $m$-shifted version. Incorporating both of these, we define
\begin{equation*}
\begin{split}
\PCar^n(Y/G,m) \, &= \, \Homcom_A( \, \Sym^n_A( \OCar^1(A)[m+1] ) \, , \, A \,)^G  \\
\, &= \, \bigoplus_{p+q = n} ( \, \Pol^p(A,m) \otimes \CoSym^q(\frg[-m]) \,)^G
\end{split}
\end{equation*}
where $\Pol^p(A,m)$ is the $m$-shifted $p$-polyvectors of $A$. 
Recall that $\Pol^p(A,m)$ has a Lie bracket of degree $-m-1$, given by the Schouten-Nijenhuis bracket. It turns out that, just as in the case of Cartan-de Rham complex, the Schouten-Nijenhuis bracket on $\Pol^p(A,m)$ induces a bracket on $\PCar^n(Y/G,m)$ so that $\PCar^*(Y/G,m)[m+1]$ becomes a weight graded DG Lie algebra (see Theorem \ref{bracket_G_inv_part}). 
One can then define an $m$-shifted Poisson structure as a Maurer-Cartan element of the DG Lie algebra $\prod_{p \geq 2} \PCar^p(Y/G,m)[m+1]$ (see Definition \ref{shifted_Poiss_def}). 
In other words, an $m$-shifted Poisson structure in Cartan model consists of an infinite sum $\pi = \pi_2 + \pi_3 + \ldots$, where $\pi_p$ is a degree $m+2$ element in $\PCar^p(Y/G,m)$, satisfying the Maurer-Cartan equation $\partial \pi + \frac{1}{2}\{\pi,\pi\} = 0$.
Just like how the above definition of a ``closed $2$-form'' relaxes the condition $d(\omega_2) = 0$, we can also think of this present definition as relaxing the condition of a Poisson bivector $\{\pi_2,\pi_2\}=0$ to the condition $\{\pi_2,\pi_2\}=-2\partial \pi_3$, together with a sequence of higher coherence relations.

Of course, what we have explained so far are only guesses. The notions of shifted symplectic and Poisson structures (especially the latter) are defined in an indirect and quite complicated way. While we think that our elementary definition is interesting as it stands, it will be more satisfying to prove that it is equivalent to the ones in the literature.
Our main result shows that this is indeed the case. More precisely, we compare with the corresponding notions in \cite{Pri17}, since it is cast in a setting close to ours.  \cite{Pri17} also suggested how to compare the definitions there with the ones in \cite{PTVV, CPTVV}. 

\bthm[=Theorems \ref{shifted_sympl_equiv}, \ref{shifted_Poiss_equiv}] \label{main_thm_intro}
For quotients of a derived affine scheme by a reductive group, the notions of shifted symplectic and Poisson structures in Cartan model coincide with the ones in \cite{Pri17}.
\ethm

As an application, we construct in \cite{Yeu} shifted symplectic and Poisson structures on derived moduli stacks of representations of (pre-)Calabi-Yau algebras using an explicit trace map.

Let us also mention that, when $G$ is not necessarily reductive, there is an extension of the Cartan-de Rham complex constructed by Getzler in \cite{Get94}, where the process of taking $G$-invariants is replaced by taking homotopy $G$-invariants (see Theorem \ref{Getzler_thm}). 
In this case, one can also define a notion of shifted symplectic structures based on the Cartan-Getzler-de Rham bicomplex (see Definition \ref{shifted_sympl_CG_def}). We also show that it is equivalent to the general notion in \cite{Pri17} (see Theorem \ref{shifted_sympl_equiv}). 

We now indicate the major technical inputs into the proof of Theorem \ref{main_thm_intro}. This will also serve as a summary of content of this paper. First, in Section \ref{Cartan_sec}, we define carefully the Cartan-de Rham bicomplex $\DRCar^*(Y/G)$ and Getzler's extension $\DR_{\CG}^*(Y/G)$, as well as the invariant $m$-shifted Cartan polyvectors $\PCar^*(Y/G,m)$. Our constructions naturally come with maps
\begin{equation}  \label{Cartan_to_CE_intro}
	\DRCar^*(Y/G) \ra \DR_{\CG}^*(Y/G) \ra \DR^*([Y/\frg])
	\qquad \text{ and } \qquad 
	\PCar^*(Y/G,m) \ra \Pol^*([Y/\frg],m)
\end{equation}
to the ordinary de Rham bicomplex and $m$-shifted polyvectors of the Chevalley-Eilenberg CDGA $\cO([Y/\frg]) := \CE(\frg^*,A)$, commuting respectively with the de Rham differential and the Schouten-Nijenhuis bracket.

The $G$-action on $Y$ induces for each $n \geq 0$ an action of $G_n := G^{n+1}$ on $X_n := Y \! \times \! G^n$, which fit into a simplicial system $[n] \in \Delta^{\op}$ (see the beginning of Section \ref{comparison_sec}). The notion of shifted symplectic and Poisson structures in \cite{Pri17} is based on the considerations of $\DR^*([X_n/\frg_n])$ and $\Pol^*([X_n/\frg_n],m)$, with a suitable notion of functoriality as we vary $[n] \in \Delta$. The proof of Theorem \ref{main_thm_intro} can be described as a $3$-step process: 
\begin{enumerate}
	\item  For each $n$, consider the map \eqref{Cartan_to_CE_intro} for the action $[X_n/G_n]$.
	\item  Show that these maps satisfy suitable functoriality as we vary $[n] \in \Delta$.
	\item  Show that it induces a quasi-isomorphism as we pass to homotopy limit over $[n] \in \Delta$.
\end{enumerate}
We have already seen Step 1. Step 2 is automatic for the de Rham side, but is non-trivial for polyvector fields. A crucial observation here is that in forming the homotopy limit over $[n] \in \Delta$, it suffices to take a subcategory $\Delta_{{\rm inj}} \subset \Delta$ which is homotopy initial (see Lemma \ref{Delta_inj_ho_init}). For a morphism $\varphi : [m] \ra [n]$ in $\Delta_{{\rm inj}}$, the corresponding map of global quotients $\varphi^* : [X_n/G_n] \ra [X_m/G_m]$ is a triangular extension in the sense of Definition \ref{triang_ext_def}. For a triangular extension, the functoriality of polyvectors in terms of a span is always well-behaved, so that we do not have to perform a fibrant replacement as in \cite{Pri17}. We provide the details of this in Section \ref{sec_functoriality}.

For Step 3, we work in a setting, established in Appendices \ref{app_cart} and \ref{app_small_alg}, where we interpret the Chevalley-Eilenberg complex as a certain ``$1$-thin quotient'' (see Proposition \ref{Gamma_star_1thin}). Accordingly, the maps \eqref{Cartan_to_CE_intro} for $[X_n/G_n]$ are interpreted as coming from a certain DG module $\OCar^p(\bX)$ over a bi-cosimplicial CDGA $\cO(\bX)$, which allows us to have a handle over the homotopy limit. Details of this are in Section \ref{comparison_sec}.

\vspace{0.2cm}

\textbf{Notations and conventions.}

For a group $G$ acting on a set $Y$ (we always consider right actions in this paper), denote by $[Y/G]_{\bullet}$ the simplicial set $[Y/G]_n = Y \times G^n$, with 
\begin{equation}  \label{Y_mod_G_def}
\begin{split}
d_0(y,g_1,\ldots,g_n) &= (yg_1,g_2,\ldots,g_n) \\
d_i(y,g_1,\ldots,g_n) &= (y,g_1,\ldots,g_{i}g_{i+1},\ldots,g_n) \, , \qquad   1 \leq i \leq n-1 \\
d_n(y,g_1,\ldots,g_n) &= (y,g_1,\ldots,g_{n-1}) \\
s_i(y,g_1,\ldots,g_n) &= (y,g_1,\ldots,g_i,e,g_{i+1},\ldots,g_n)  \, , \qquad  0 \leq i \leq n
\end{split}
\end{equation}
In other words, $[Y/G]_{\bullet}$ is the nerve of the action groupoid of the left action of $G^{\op}$ on $Y$.

We fix a field $k$ of characteristic zero. Unadorned tensor products are understood to be over $k$. We work with cochain complexes (\ie, differentials have degree $1$). Shifts are defined by $M[1]^n = M^{n+1}$. The category of cochain complexes of $k$ vector spaces is denoted by $\Ch(k)$.
For DG modules $M,N$ over a DG algebra $A$, their Hom complex will be denoted by $\Homcom_A(M,N)$. Hence $\Hom_A(M,N) = Z^0(\Homcom_A(M,N))$.

The notations $\DGA_k$, $\CDGA_k$, $\CDGAn_k$, $\CDGAp_k$, $\dgMod(A)$, $\dgMod_G(A)$, $\Alg_k$, $\CAlg_k$, $\setdel$ will refer respectively to the categories of differential graded algebras (DGA), commutative DG algebras (CDGA), non-positively graded CDGA, non-negatively graded CDGA, DG modules over $A$, $G$-equivariant DG modules over $A$, associative algebras, commutative algebras, and simplicial sets.
The category $\dAff$ of derived affine schemes is defined to be the category opposite to $\CDGAn_k$.
The notation $\cC^I$ will refer to the category of functors from $I$ to $\cC$. For example, $\CDGA^{\Delta}_k$ means the category of cosimplicial CDGAs.

The categories $\Ch(k)$, $\DGA_k$, $\CDGA_k$ and $\dgMod(A)$ (resp.  $\Ch^{\leq 0}(k)$ and $\CDGAn_k$) are understood to be endowed with the standard model structures, where a map $f$ is a weak equivalence if and only if it is a quasi-isomorphism, and is a fibration if and only if it is surjective in all degrees (resp. in positive degrees). Given a simplicial object $M_{\bullet}$ (resp. cosimplicial object $M^{\bullet}$) in $\Ch(k)$ or $\dgMod(A)$, we will denote by $\Tot^{\oplus}_{\Delta^{\op}} \, M_{\bullet}$ (resp. $\Tot^{\Pi}_{\Delta} \, M^{\bullet}$) the direct sum (resp. direct product) total complex of its associated unnormalized bicomplex. We will refer to $\Tot^{\oplus}_{\Delta^{\op}} \, M_{\bullet}$ as the homotopy colimit and $\Tot^{\Pi}_{\Delta} \, M^{\bullet}$ as the homotopy limit, as it is an explicit model for such.

We also specify our convention for the Chevalley-Eilenberg complex here.
We always write $\Sym^n_k(V)$ for the $S_n$-coinvariants on $V^{\otimes n}$ under (Koszul signed) permutation. Unless explicitly specified, we will not identify it with the $S_n$-invariants. The structure map $\delta_{\frg^*} : \frg^* \ra \frg^* \otimes \frg^*$ of a Lie coalgebra lands in the ($S_2$-invariant part of) $\frg^* \otimes \frg^*$. The notion of a $\frg^*$-comodule $M$ dualizes that of a Lie algebra module $N$, defined by $\cdot : N \otimes \frg \ra N$ satisfying $(\xi \cdot x) \cdot y  - ( \xi \cdot y) \cdot x =  \xi \cdot [x,y]$ (we do not include a factor $\tfrac{1}{2}$ on the left hand side). Given a $\frg^*$-comodule $M$ with structure map $\beta_M : M \ra  M \otimes \frg^*$, then the Chevalley-Eilenberg complex is given by $\CE(\frg^*,M) = M \otimes \Sym(\frg^*[-1])$ with differential $\delta$ given by $\delta|_{\frg^*[-1]} : \frg^*\xra{\frac{1}{2}\delta_{\frg^*}} \frg^* \otimes \frg^* \ronto \Lambda^2 \frg^*$ and $\delta|_M = \delta_M$. 
More precisely, the differential is determined by the above specification, as well as by the requirement that $\CE(\frg^*,k)$ is a CDGA, and $\CE(\frg^*,M)$ is a DG module over $\CE(\frg^*,k)$.
Notice that the factor $\tfrac{1}{2}$ on $\delta|_{\frg^*[-1]}$ is necessary to ensure that $\delta^2 = 0$ on $\CE(\frg^*,M)$.
In the case when $A$ is a commutative algebra with a $\frg^*$-comodule algebra structure, $\CE(\frg^*,A)$ is a CDGA, which we will write as $\cO([Y/\frg])$. 

All of the above discussions for the Chevalley-Eilenberg complex have a DG version. For example, if $A$ is a $\frg^*$-comodule CDGA, then $\cO([Y/\frg])$ is a bigraded CDGA. Similarly, in Appendix \ref{app_small_alg}, results will be stated only for commutative algebras, but will be applied to CDGAs in the main text without further mention.

\section{Shifted symplectic and Poisson structures in Cartan and Cartan-Getzler models}  \label{Cartan_sec}

Let $G$ be an linear algebraic group (automatically smooth since we assume ${\rm char}(k)=0$) acting on a non-positively graded CDGA $(A,\partial)$, with coaction map $\Delta : A \ra A \otimes \cO(G)$, then we may take the composition
\begin{equation*}
\beta \,: \, A \xraq{\Delta}  A \otimes \cO(G) \xraq{\id \otimes d_e}   A \otimes \frg^*
\end{equation*}
where $d_e : \cO(G) \ra \frg^*$ sends $f$ to the differential of $f$ at $e \in G(k)$. 

This map $\beta$ is a derivation on the free $A$-module $A\otimes \frg^*$, and hence descends to an $A$-linear map
\begin{equation}  \label{alpha_Omega1}
\alpha \, : \, \Omega^1(A) \raq A \otimes \frg^*
\end{equation}

We will take the fiber (or cocone) of this map (with suitable gradings to be specified presently)
\begin{equation}  \label{OCar1_cone}
\OCar^1(\frg^*,A) \, := \, [\, \Omega^1(A) \xraq{\alpha}   A \otimes \widetilde{\frg^*} \, ]
\end{equation}
where we write $\widetilde{\frg^*}$ for $\frg^*$ in order to avoid possible notational confusion below.
Here, we introduce two new gradings: one is the equivariant grading, and one is the Hodge grading. The fiber is to be taken in the equivariant grading, so that $\OCar^1(\frg^*,A)$ is a bicomplex with $\Omega^1(A)$ sitting in equivariant degree $0$ and $A \otimes \widetilde{\frg^*}$ in equivariant degree $1$, both of which inherits intrinsic gradings and differentials from $(A,\partial)$. We put $\OCar^1(\frg^*,A)$ in Hodge degree $1$.

Then we take
\begin{equation*}
\OCar^p(\frg^*,A) \, := \, \Sym_A^p(\, \OCar^1(\frg^*,A) \,)
\end{equation*}
Notice that the Koszul rule for multiply graded complexes is given by $\sigma(x\otimes y) = (-1)^{\sum_{i=1}^r m_i n_i} y \otimes x$, for multi-gradings $|x|= (m_1,\ldots,m_r)$ and $|y|= (n_1,\ldots,n_r)$.
If desired, one can always collapse several of the gradings by taking the total complex.
The collapsing operation is strong symmetric monoidal. For example, the functor $\Tot$ that collapse a double complex to a single complex can be enriched to a strong monoidal functor by the maps $\Tot(M \otimes N) \xra{\cong} \Tot(M) \otimes \Tot(N)$ given by $x \otimes y \mapsto (-1)^{qr} x\otimes y$ for $|x|=(p,q)$ and $|y| = (r,s)$.

Explicitly, $\OCar^*(\frg^*,A)$ is given by 
\begin{equation}  \label{OCar_Omega_Sym}
\OCar^*(\frg^*,A) \, \cong \,   \Omega^*(A) \otimes \Sym^*(\widetilde{\frg^*}[0,-1,-1])
\end{equation}
where we have written the shift in the tri-grading in the order (intrinsic, equivariant, Hodge), and $\Omega^p(A)$ is understood to lie in Hodge degree $p$.

Both of the $A$-modules in \eqref{OCar1_cone} defining the cone have canonical $G$-equivariant structures induced from $A$ and $\frg^*$. The map $\alpha$ between them is also $G$-equivariant, so that $\OCar^1(\frg^*,A)$ has a $G$-equivariant structure. In other words, the tri-graded CDGA $\OCar^*(\frg^*,A)$ has a $G$-action extending the given one on $A$.

\bdf
The bigraded complex $\OCar^p(\frg^*,A)$ is called the \emph{bicomplex of Cartan $p$-forms}.
The $G$-invariant part $\OCar^p(\frg^*,A)^G$ is called the \emph{bicomplex of invariant Cartan $p$-forms}.
\edf

Denote the de Rham differential on $\Omega^*(A)$ by $d$, and consider the differential $d' = d \otimes \id$ on the right hand side of \eqref{OCar_Omega_Sym}. The differential $d'$ clearly commutes with the intrinsic differential $\partial$, but in general does not commute with the equivariant differential (induced from the cone \eqref{OCar1_cone}), which we denote as $\delta'$. However, it is a classical fact that it does commute with $\delta'$ on the $G$-invariant part:
\bthm  \label{Omega_CE_dR_inv}
The differential $d'$ commutes with $\delta'$ when restricted to the invariant part $\OCar^*(\frg^*,A)^G$.
\ethm

Theorem \ref{Omega_CE_dR_inv} is usually verified directly. Alternatively, it also follows from a comparison with a certain Chevalley-Eilenberg complex. First, since $\OCar^*(\frg^*,A)$ has a $G$-action, it has an induced $\frg^*$-comodule algebra structure, so that we may form the Chevalley-Eilenberg complex $\CE(\frg^*,\OCar^*(\frg^*,A))$. This is then a quadri-graded complex, with a new grading and differential, which will be called the second equivariant grading and differential $\delta''$. As a quadri-graded commutative algebra, it is given by 
\begin{equation}  \label{CE_OCar}
\CE(\frg^* , \OCar^*(\frg^*,A)) \, \cong \,     \Omega^*(A) \, \otimes \,  \Sym(\frg^*[0,(0,-1),0]) \, \otimes \, \Sym(\widetilde{\frg^*}[0,(-1,0),-1])
\end{equation}
where we have written the gradings according to (intrinsic, (first equivariant, second equivariant), Hodge).

On the other hand, we can take the Chevalley-Eilenberg complex of $A$ itself, which we will denote as $\cO([Y/\frg]) = \CE(\frg^*,A)$, where $Y = \Spec \, A$. This is a bi-graded CDGA, with (intrinsic, equivariant) gradings. Then, we take its de Rham complex $\Omega^*([Y/\frg])$, which is a tri-graded CDGA with (intrinsic, equivariant, Hodge) gradings. As a tri-graded commutative algebra, it is given by 
\begin{equation}  \label{Omega_Yg}
\Omega^*([Y/\frg]) \, = \,     \Omega^*(A)  \otimes \Sym(\frg^*[0,-1,0])  \otimes \Sym(\widetilde{\frg^*}[0,-1,-1])
\end{equation}
where the de Rham differential $d$ sends the copy $\frg^*$ to the copy $\widetilde{\frg^*}$ by the identity map, and is the usual one on $\Omega^*(A)$.

Comparing the graded commutative algebras \eqref{CE_OCar} and \eqref{Omega_Yg}, we see that they are the same if we collapse the two equivariant gradings in \eqref{CE_OCar}. The intrinsic differential also clearly match in this identification. In fact, the equivariant differentials also match:
\bpp  \label{CE_Omega_Car_Omega_CE}
If we collapse the first and second equivariant gradings in $\CE(\frg^* , \OCar^*(\frg^*,A))$, then the identification $\CE(\frg^* , \OCar^*(\frg^*,A)) \cong \Omega^*([Y/\frg])$ given by \eqref{CE_OCar} and \eqref{Omega_Yg} commutes with the equivariant differentials.
\epp

\bpf
The $G$-action on $A$ and the induced $G$-action on $\Omega^*(A)$ induces $\frg^*$-comodule algebra structures $\beta_A : A \ra A \otimes \frg^* $ and $\beta_{\Omega^*} : \Omega^*(A) \ra \Omega^*(A) \otimes \frg^*$ by post-composing the coaction maps with $d_e : \cO(G) \ra \frg^*$. The equivariant differential $\delta$ on $\Omega^*([Y/\frg])$ can be described in terms of these maps. Namely, using the convention we specified at the end of Section \ref{sec_intro},
one can verify directly
that the restriction of $\delta$ on $A$, $\Omega^1(A)$, $\frg^*$ and $\widetilde{\frg^*}$ are given by 
\begin{equation}  \label{delta_Omega_CE}
\begin{split}
\delta|_A \, &= \, \beta_A \, : \, A \raq A \otimes \frg^* \\
\delta|_{\Omega^1(A)} \, &= \,  ( \,\beta_{\Omega^1} \, , \, \alpha \,) \, : \, \Omega^1(A) \raq ( \Omega^1(A) \otimes \frg^* ) \,  \oplus \,(A \otimes \widetilde{\frg^*}) \\
\delta|_{\frg^*} \, &: \, \frg^*\xra{\frac{1}{2}\delta_{\frg^*}} \frg^* \otimes \frg^* \ronto \Lambda^2 \frg^* \\
\delta|_{\widetilde{\frg^*}} \, &: \, \widetilde{ \frg^*} \,=\, \frg^* \xraq{\delta_{\frg^*}}  \frg^* \otimes  \frg^*  \,= \,  \frg^* \otimes \widetilde{ \frg^*}
\end{split}
\end{equation}
where $\delta_{\frg^*}$ is the Lie coalgebra structure map.
This clearly matches with $\delta' + \delta''$ on $\CE(\frg^* , \OCar^*(\frg^*,A))$.
\epf

\bpf[Proof of Theorem \ref{Omega_CE_dR_inv}]
Translate the de Rham differential on $\Omega^*(Y/\frg)$ to a de Rham differential $d$ on $\CE(\frg^* , \OCar^*(\frg^*,A))$. 
It is easy to see that $d = d' + d''$ where $d' = d \otimes \id \otimes \id $ under the presentation \eqref{CE_OCar}, and $d''$ is the unique derivation that sends $\frg^*$ to $\widetilde{ \frg^*}$ via the identity map, and is zero on $\Omega^*(A)$ and $\widetilde{ \frg^*}$. By Proposition \ref{CE_Omega_Car_Omega_CE}, we see that $d' + d''$ commutes with $\delta' + \delta''$. 
Notice that, with respect to the first and second equivariant gradings, $d'$ has degree $(0,0)$ while $d''$ has degree $(1,-1)$.
Since $d'+d''$ commutes with $\delta' + \delta''$, we see in particular that $d'$ commutes with $\delta''$ for degree reasons.  
Denote by $C^0$ and $Z^0$ the degree $0$ elements and cocycles with respect to the second equivariant differential $\delta''$. \ie, $C^0 = \OCar^*(\frg^*,A)$ and $Z^0 = \OCar^*(\frg^*,A)^{\frg}$. Then by degree reasons $d'$ preserves $C^0$, and hence $Z^0$ since we have just seen that $[d',\delta'']=0$.
Since $d''=0$ and $\delta'' = 0$ on $Z^0$, we see that $d'$ commutes with $\delta'$ on $Z^0 = \OCar^*(\frg^*,A)^{\frg}$, and hence in particular on $\OCar^*(\frg^*,A)^{G}$.
\epf

Thus, the invariant Cartan forms $\OCar^*(\frg^*,A)^G$ becomes a tri-complex, with differentials denoted as $(\partial, \delta',d')$, along the (intrinsic, equivariant, Hodge) gradings respectively.
We will rewrite it as $\DRCar^*(\frg^*,A) := \OCar^*(\frg^*,A)^{G}$

\bdf
The tri-complex $(\DRCar^*(\frg^*,A) ,\partial, \delta', d')$ is called the \emph{Cartan-de Rham tri-complex}.
\edf

Next we consider the dual picture of polyvector fields. 
For any $B \in \CDGA_k$, let the complex of $m$-shifted $p$-polyvector fields be
\begin{equation*}
\Pol^p(B,m) \, := \, \Homcom_B( \, \Sym^p_B(\Omega^1(B)[m+1]) \, , \, B \,)
\end{equation*}
Notice that, unlike in the case for differential forms, we do not introduce a separate Hodge grading (alternatively, it means that we have already collapsed the Hodge degree with the intrinsic one from the start). At times, we may say that $\Pol^p$ has weight grading $p$, but such a weight grading does not contribute to Koszul signs.

We recall that $\Pol^*(B,m)$ has a canonical DG $P_{m+2}$-algebra structure (see also \cite{CPTVV}).
The product structure is a shuffle product. Namely, if $M$ is any DG module over $B$, then $\Sym_B^*(M)$ has a $B$-linear bialgebra structure where the product is the usual one, and the (cocommutative) coproduct is determined by $\Delta(\xi) = \xi \otimes 1  + 1 \otimes \xi$ for $\xi \in M$.
Taking the $B$-linear dual, $\Homcom_B(\Sym_B^*(M),B)$ then inherits a commutative product (which is given by a sum over shuffle permutations).

Explicitly, we write $\langle n \rangle := \{1,\ldots,n\}$, regarded as a totally ordered set. For any $\theta \in M^{\otimes p+q}$, and any partition $\langle p+q \rangle = S' \amalg S''$ for $|S'| = p$ and $|S''| = q$, there exists a unique $\sigma \in S_{p+q}$ that is order-preserving on both $S'$ and $S''$, and with images $\sigma(S') = \{1,\ldots,p\}$ and $\sigma(S'') = \{p+1,\ldots,p+q\}$.
Write in Sweedler's notation $\sigma(\theta) = \theta_{S'} \otimes \theta_{S''} \in M^{\otimes p} \otimes M^{\otimes q}$.
Then the product on $\Homcom_B(\Sym_B^*(M),B)$ is given by
\begin{equation}  \label{shuffle_prod_Hoch}
(f \cdot g)(\theta) \, = \, \sum_{\substack{\langle p+q \rangle = S' \amalg S'' \\ |S'| = p \text{ and }|S''| = q}}	f(\theta_{S'}) \cdot g(\theta_{S''})
\end{equation}

The Schouten-Nijenhuis bracket is defined by considering $\Pol^*(B,m)$ as a subcomplex
\begin{equation*}
\Homcom_B( \, \Sym^p_B(\Omega^1(B)[m+1]) \, , \, B \,) \, \subset \, \Homcom_k(  (B[m+1])^{\otimes p}  , B )
\end{equation*}
consisting of those maps that are a derivation in each variable, and is symmetric under the $S_p$-action.
The right hand side (shifted by $m+1$) has a dg pre-Lie algebra structure given by
\begin{equation}  \label{Hoch_pre_Lie_1}
	f * g \, := \, \sum_{i=1}^p \, \sum_{\sigma \in S^i_{p,q}} \,(f \circ_i g)^{\sigma}
\end{equation}
for $f \in \Homcom_k(  (B[m+1])^{\otimes p}  , B[m+1] )$ and $g \in \Homcom_k(  (B[m+1])^{\otimes q}  , B[m+1] )$, where 
\begin{equation*}
	S^i_{p,q} \,:= \, \Bigl\{ \, \sigma \in S_{p+q-1} \, | \, \substack{\sigma^{-1}(1)<\sigma^{-1}(2) <\ldots< \sigma^{-1}(i) < \sigma^{-1}(i+q) < \ldots < \sigma^{-1}(p+q-1) \\
		\text{and } \sigma^{-1}(i) < \sigma^{-1}(i+1) < \ldots < \sigma^{-1}(i+q-1)} \, \Bigr\}
\end{equation*}



For any subset $S \subset \langle p+q-1 \rangle$ of cardinality $|S|=q$, 
there exists a unique pair $(i,\sigma)$ where $1 \leq i \leq p$ and $\sigma \in S^i_{p,q}$, such that $\sigma(S) = \{i, i+1,\ldots,i+q-1\}$. (Here, $i$ is simply the smallest element of $S$). Conversely, the pair $(i,\sigma)$ uniquely determines $S = \sigma^{-1}(\{i, i+1,\ldots,i+q-1\})$.
For any $1 \leq j \leq p$, denote by $g_{S,j} : (B[m+1])^{\otimes p+q-1} \ra (B[m+1])^{\otimes p}$ the map that sends $(B[m+1])^{\otimes S}$ to the $j$-th copy of $B[m+1]$ via $g$ ($S$ inherits a total order from $\langle p+q-1 \rangle$, so that $(B[m+1])^{\otimes S}$ is canonically identified with $(B[m+1])^{\otimes q}$), and the identity map on the rest (in an order-preserving way on the tensor indexing). Write $g_S := g_{S,j}$ for $j = \min(S)$.
Then \eqref{Hoch_pre_Lie_1} can be rewritten as
\begin{equation} \label{Hoch_pre_Lie_2}
	f * g \, := \, \sum_{S \subset \langle p+q-1 \rangle, \, \, |S|=q} f \circ (g_{S})
\end{equation}

One can show that the associated Lie bracket $\{f,g\} := f * g - (-1)^{|f||g|} g * f$ preserves the subcomplex $\Pol^*(B,m)[m+1]$ (see, {\it e.g.}, \cite[Section 2]{Mel16}), which we take as the Schouten-Nijenhuis bracket $\{-,-\}$. One can verify that it is a derivation in each variable with respect to the shuffle product \eqref{shuffle_prod_Hoch}, so that it makes $\Pol^p(B,m)$ into a $P_{m+2}$-algebra.
The Schouten-Nijenhuis bracket is the usual commutator bracket on weights $1$, and is the usual action of weight $1$ on weight $0$, so that there is always a map of $P_{m+2}$-algebras $\Sym_B(\Der(B)[-m-1]) \ra \Pol^*(B,m)$, which is moreover an isomorphism if $\Omega^1(B)$ is finitely generated and projective as a graded module over $B$.

Now, consider again a linear algebraic group $G$ acting on a non-positively graded CDGA $(A,\partial)$. Let 
\begin{equation*}
\XCar^p(\frg, A,m) \, := \, \Homcom_A( \, \Sym^p_A(\OCar^1(\frg^*,A)[m+1]) \, , \, A \,)
\end{equation*}
which is a bicomplex (with bounded non-positive equivariant degrees). In other words,
\begin{equation}  \label{PCar_Pol_g}
\XCar^*(\frg,A,m) \, \cong \,  \Pol^*(A,m) \otimes \CoSym^*(\widetilde{\frg}[-m-1,1]) 
\end{equation}
where we have written $\widetilde{\frg}$ for $\frg$ to avoid possible notational confusion below. We have written the shifts in (intrinsic, equivariant) gradings. 
Explicitly, the identification \eqref{PCar_Pol_g} is obtained by a map from the right to the left, given as a shuffle sum.

Notice that if $M,N$ are $G$-equivariant modules, then $\Homcom_A(M,N)$ inherits a $G$-equivariant module structure only under certain finiteness assumption (for example, it holds if $M$ is projective of finite rank as a graded module). However, in general, it still has enough remnants of a $G$-equivariant structure for all our purposes. Firstly, it still has a $\frg^*$-comodule structure (we assume that $G$ is a linear algebraic group, so that $\frg^*$ is finite dimensional), defined by
\begin{equation}  \label{Hom_g_comod}
\beta : \Homcom_A(M,N) \ra  \Homcom_A(M,N) \otimes \frg^* \, , \qquad \beta(f) = \beta_M \circ f - (f \otimes \id) \circ \beta_M
\end{equation}
so that we may define its Chevalley-Eillenberg complex $\CE(\frg^*, \Homcom_A(M,N))$. Secondly, one can still define a ``complex of homotopy $G$-invariants'', which is a module $C^{\bullet}(G,\Homcom_A(M,N))$ over the cosimplicial CDGA $\cO([Y/G]_{\bullet})$, given by $C^n(G,\Homcom_A(M,N)) = \Homcom_A(M,N \otimes \cO(G)^{\otimes n})$ (see \eqref{Hom_sim}).
In particular, $\Homcom_A(M,N)^G = \ker(d^0-d^1)$ is well-defined, and consists of the $G$-equivariant maps from $M$ to $N$. 
Thirdly, if we think of $\Homcom_A(M,N)$ as a functor $\CAlg_k \ra \Ch(k)$ given by $R \mapsto \Homcom_A(M,N \otimes R)$, then this functor has an action by the group valued functor $R \mapsto G(R)$.
We may apply these considerations to $\XCar^p(\frg, A,m)$, which then inherits a ``$G$-equivariant structure'' in any of these senses.

Define the bracket $\{-,-\}'$ on $\XCar^*(\frg,A,m)$ by $\{-,-\}' = \{-,-\} \otimes \mu$ on the right hand side of \eqref{PCar_Pol_g}, where $\mu$ is the (shuffle) product map. It is clearly a graded Poisson bracket of (intrinsic) degree $-m-1$ and equivariant degree $0$. It clearly commutes with the intrinsic differential, but it does not commute with the equivariant differential in general. However, as in the case for Cartan forms, we have the following
\bthm  \label{bracket_G_inv_part}
The bracket $\{-,-\}'$ commutes with the equivariant differential on the $G$-invariant part $\XCar^*(\frg,A,m)^G$, making it into a bigraded DG $P_{m+2}$-algebra (more precisely, the bracket has degree $(-m-1,0)$ with respect to (intrinsic, equivariant) grading).
\ethm

As in the case for Cartan forms, we will prove this by comparing with a certain Chevalley-Eilenberg complex. First, note that
\begin{equation}  \label{CE_PCar}
\CE(\frg^*, \XCar^*(\frg,A,m)) \, \cong \,  \Pol^*(A,m) \otimes  \Sym^*(\frg^*[0,(0,-1)]) \otimes  \CoSym^*(\widetilde{\frg}[-m-1,(1,0)])
\end{equation}
where we have introduced a new grading (the second equivariant grading), and we now writes the gradings in the order (intrinsic, (first equivariant, second equivariant)).

On the other hand, we can also form $\Pol([Y/\frg],m)$, which can be described as
\begin{equation}  \label{PCar_Y_g}
\Pol^*([Y/\frg],m) \, \cong \,   \Pol^*(A,m) \otimes \Sym^*(\frg^*[0,-1]) \otimes  \CoSym^*(\widetilde{\frg}[-m-1,1])
\end{equation}

We have the following analogue of Proposition \ref{CE_Omega_Car_Omega_CE}:

\bpp  \label{Pol_CE_isom}
If we collapse the first and second equivariant gradings in $\CE(\frg^*, \XCar^*(\frg,A,m))$, then the identification $\CE(\frg^*, \XCar^*(\frg,A,m)) \cong \Pol^*([Y/\frg],m)$ given by \eqref{CE_PCar} and \eqref{PCar_Y_g} commutes with the equivariant differentials.
\epp

\bpf
This follows from a direct computation. For $\Pol^*([Y/\frg],m)$, take the $\cO([Y/\frg])$-linear dual of \eqref{delta_Omega_CE}. For $\CE(\frg^*, \XCar^*(\frg,A,m))$, use \eqref{Hom_g_comod} to obtain the $\frg^*$-comodule structure.
\epf

\bpf[Proof of Theorem \ref{bracket_G_inv_part}] 
Transport the Schouten-Nijenhuis bracket on $\Pol^*([Y/\frg],m)$ to $\CE(\frg^*, \XCar^*(\frg,A,m))$. Under the presentation \eqref{CE_PCar}, the translated bracket is the sum of two parts $\{-,-\}' + \{-,-\}''$, given by
$\{-,-\}' \, = \, \{-,-\} \otimes \mu $ and $
\{-,-\}'' \, = \, \mu \otimes \{-,-\}_{\frg} $,
where $\{-,-\}$ is the Schouten-Nijenhuis bracket on $\Pol^*(A,m)$, and $ \{-,-\}_{\frg}$ is the derivation that extends the linear pairing on $\frg$ and $\frg^*$, and $\mu$ are the multiplication maps on the remaining parts.

By Proposition \ref{Pol_CE_isom}, $\{-,-\}' + \{-,-\}''$ commute with $\delta' + \delta''$.
Notice that $\{-,-\}'$ and $\{-,-\}''$ have equivariant degrees $(0,0)$ and $(1,-1)$ respectively. The rest of the proof then goes exactly as in the proof of Theorem \ref{Omega_CE_dR_inv}.
\epf

Now we will use these to define shifted symplectic and Poisson structures. So far, for bookkeeping purposes, we have been careful to separate the intrinsic, first, and second equivariant gradings. Once all the above work has been done, we will now collapse all of them into the intrinsic degree and differential.

\bdf
The complex $\OCar^p(Y/G) := \Tot \, \OCar^p(\frg^*,A)$ is called the \emph{complex of Cartan $p$-forms}.
The $G$-invariant part $\OCar^p(Y/G)^G = \Tot \, \OCar^p(\frg^*,A)^G$ is called the \emph{complex of invariant Cartan $p$-forms}, which we rewrite as 
$\DRCar^p(Y/G) := \OCar^p(Y/G)^G$.
The bicomplex $(\DRCar^*(Y/G) ,\partial, d')$ is called the \emph{Cartan-de Rham bicomplex}.

The complex $\XCar^p(Y/G,m) := \Tot \, \XCar^p(\frg, A,m)$ is called the \emph{complex of $m$-shifted Cartan $p$-polyvectors}. 
The $G$-invariant part $\XCar^p(Y/G,m)^G := \Tot \, \XCar^p(\frg, A,m)^G$ is called the \emph{complex of $m$-shifted invariant Cartan $p$-polyvectors}, which we rewrite as 
$\PCar^p(Y/G,m) = \XCar^p(Y/G,m)^G$. 
The Lie bracket $\{-,-\}'$ on $\PCar^*(Y/G,m)$ is called the \emph{Schouten-Nijenhuis-Cartan bracket}, which makes it into a (weight graded) DG $P_{m+2}$-algebra.

Notice that in all these formations of total complexes, since we fix a Hodge degree, there are only finitely many equivariant degrees with non-zero terms. Hence the direct product total complex coincides with the direct sum total complex, so there is no ambiguity.
\edf

We may now proceed to define shifted symplectic structures. As explained in Section \ref{sec_intro}, in order to define closed forms, we further collapse the Cartan-de Rham bicomplex by taking a direct product total complex
\begin{equation*}
\DRCarcl(Y/G) \, := \, \Pi_{p \geq 0} \, \, \DRCar^p(Y/G)[-p]
\end{equation*}
which comes with the Hodge filtration
\begin{equation*}
F^n \, \DRCarcl(Y/G) \, := \, \Pi_{p \geq n} \, \, \DRCar^p(Y/G)[-p]
\end{equation*}

It is clear that any cocycle $\omega \in Z^m(\DR^2(Y/G))$ induces a 
$G$-equivariant%
\footnote{As we mentioned above in the paragraph of \eqref{Hom_g_comod}, $\OCar^1(Y/G)^{\vee}$ has a $G$-equivariance structure in several senses. This map is $G$-equivariant in any of these senses.} map
\begin{equation}  \label{omega_sharp}
\omega^{\sharp} \, : \, \OCar^1(Y/G)^{\vee}[-m] \raq \OCar^1(Y/G)
\end{equation}
of DG modules over $A$.

\bdf  \label{almost_cof_def}
We say that $A \in \CDGAn_k$ is \emph{almost cofibrant} (or \emph{smooth}) if $\Omega^1(A) \in \dgMod(A)$ is cofibrant, and if the map $\bL_{A/k} \ra \Omega^1(A)$ is a quasi-isomorphism. 
\edf

\bdf  \label{shifted_sympl_def}
Suppose that $G$ is reductive and $A$ is almost cofibrant, then an \emph{$m$-shifted pre-symplectic structure in Cartan model} on $[Y/G]$ is a cocycle $\widetilde{\omega} \in Z^{m+2}F^2 \DRCarcl(Y/G)$.

Denote by $\omega \in Z^m (\DR^2(Y/G))$ the part of $\widetilde{\omega}$ in Hodge degree $2$ (\ie, it is the image under the quotient by $F^3$). Then $\widetilde{\omega}$ is said to be \emph{non-degenerate} if the corresponding map \eqref{omega_sharp} is a quasi-isomorphism. 
In this case, we say that $\widetilde{\omega}$ is an \emph{$m$-shifted symplectic structure in Cartan model} on $[Y/G]$.

Two $m$-shifted pre-symplectic structures in Cartan model are said to be \emph{equivalent} if they lie in the same cohomology class $H^{m+2}F^2 \DRCarcl(Y/G)$. Clearly, non-degeneracy depends only on the equivalence class.
The set of equivalence classes of (pre-)symplectic structures is in a natural bijection to the set of connected components of a certain space.
Namely, define the \emph{space of $m$-shifted pre-symplectic structures in Cartan model} to be the Dold-Kan denormalization of the good truncation $\tau^{\leq 0}( F^2 \DRCarcl(Y/G)[m+2] )$. 
Define the \emph{space of $m$-shifted symplectic structures in Cartan model} to be the union of the connected components that consist of points that are non-degenerate.
\edf

To define shifted Poisson structures, we likewise form a direct product
\begin{equation*}
\hPCar(Y/G,m) \, := \, \prod_{p \geq 0} \, \PCar^p(Y/G,m)
\end{equation*}
which comes with the weight filtration
\begin{equation*}
F^n\hPCar(Y/G,m) \, := \, \prod_{p \geq n} \, \PCar^p(Y/G,m)
\end{equation*}

The Schouten-Nijenhuis-Cartan bracket makes $\hPCar(Y/G,m)[m+1]$ a DG Lie algebra, and each of the subcomplexes $F^n\hPCar(Y/G,m)[m+1]$ a DG Lie subalgebra.

\bdf  \label{shifted_Poiss_def}
Suppose that $G$ is reductive and $A$ is almost cofibrant, then an \emph{$m$-shifted Poisson structure in Cartan model} on $[Y/G]$ is a Maurer-Cartan element in the DG Lie algebra $F^2\hPCar(Y/G,m)[m+1]$.

The \emph{space of $m$-shifted Poisson structure in Cartan model} is the Maurer-Cartan space, defined either as the mapping space $\Map_{\DGLA^{\gr}}( k[-1](-2), \PCar(Y/G,m)[m+1] )$ in the model category $\DGLA^{\gr}$ of weight graded DG Lie algebras (where the bracket has weight $-1$), or more concretely as in \cite[Definition 1.5]{Pri17}.
\edf

We now discuss shifted symplectic structures in the case when $G$ is not necessarily reductive. In this case, it is unsatisfactory to take the $G$-invariants $\OCar^p(\frg^*,A)^G$. Instead, one should take its homotopy $G$-invariants. Namely, we consider the associated Cartesian module $C^{\bullet}(G,\OCar^p(\frg^*,A))$ (see \eqref{assoc_cart}), and take its normalized complex, which we denote by $N^{\bullet}(G,\OCar^p(\frg^*,A))$.
The cosimplicial degree then gives a new grading, which we call the second equivariant grading, so that $N^{\bullet}(G,\OCar^p(\frg^*,A))$ is a tri-complex.
We again put $N^{\bullet}(G,\OCar^p(\frg^*,A))$ in Hodge grading $p$.

Recall from Proposition \ref{glob_quot_1thin_quot} that for any CDGA $B$ with a $G$-action, there is a canonical map
\begin{equation*}
\pi : N(\cO([\Spec \, B / G])) \ronto \CE(\frg^*,B)
\end{equation*}
given as the $1$-thin quotient, explicitly described in \eqref{glob_quot_1thin_quot_map}. 

In particular, if we apply this to $B = \OCar^*(\frg^*,A)$, then there is a map of tri-complexes
\begin{equation}  \label{OCar_hG_to_hg}
\pi : N^{\bullet}(G,\OCar^p(\frg^*,A))  \rontoq \CE(\frg^*,\OCar^p(\frg^*,A))
\end{equation}

We may use this as a hint to construct a de Rham differential on $N^{\bullet}(G,\OCar^p(\frg^*,A))$. Recall that we constructed above a de Rham differential on $\CE(\frg^*,\OCar^*(\frg^*,A))$ of the form $d = d' + d''$, where $d'$ has equivariant degrees $(0,0)$ and $d''$ has equivariant degrees $(1,-1)$. Moreover, $d = d' + d''$ commutes with $\delta = \delta' + \delta''$, but not independently.
One may ask whether $d'$ and $d''$ lifts to $N^{\bullet}(G,\OCar^p(\frg^*,A))$, with similar properties. By a construction of \cite{Get94}, the answer turns out to be affirmative. But for now, let's assume that we do not know the answer, and try to guess the formula. 

We will focus our attention to differentials that satisfy the following property:
\begin{equation}  \label{dd_prop}
\parbox{40em}{The maps $d'$ and $d''$ are derivations of $N^{\bullet}(G,\OCar^p(\frg^*,A))$ with respect to the cup product, and have equivariant degrees $(0,0)$ and $(1,-1)$ respectively.}
\end{equation}

If we assume this, then it turns out that the most straightforward guess works. For $d'$, we simply take
\begin{equation}  \label{dd_eqn1}
d' = d \otimes \id \otimes \id  \quad \text{ on } \quad N^n( G,\OCar^*(\frg^*,A) ) =  \Omega^*(A) \otimes \Sym(\widetilde{\frg^*}[0,(-1,0),-1]) \otimes \frm^{\otimes n}  
\end{equation}

To define $d''$, notice that $N^{\bullet}(G,B)$ is generated by $N^0(G,B)=B$ and $1 \otimes \frm \subset N^1(G,B)$ under the cup product (see Appendix \ref{app_small_alg}). Since $d''$ has equivariant degree $(1,-1)$, it must vanish on $N^0$, so that it suffices to specify it on $1 \otimes \frm$. The most straightforward guess is then given by
\begin{equation}  \label{dd_eqn2}
d''(1 \otimes f) = \overline{f} \in \widetilde{\frg^*} \quad \text{ for } \quad 1 \otimes f \,  \in \, 1 \otimes \frm \subset N^1(G,\OCar^*(\frg^*,A))
\end{equation}
where $f \mapsto \overline{f}$ is the map $ \frm \ronto \frm/\frm^2 = \frg^*$ to the quotient.

It turns out that these choices work:

\bthm[\cite{Get94}]  \label{Getzler_thm}
There are unique differentials $d'$ and $d''$ satisfying \eqref{dd_prop}, \eqref{dd_eqn1} and \eqref{dd_eqn2}. Moreover, $d = d' + d''$ commutes with $\delta = \delta' + \delta''$ and satisfies $d^2 = 0$. They lift $d'$ and $d''$ respectively under the map \eqref{OCar_hG_to_hg}.
\ethm

\bpf
See \cite[Corollary 2.1.2, Proposition 2.3.1]{Get94}. One can also simplify the proof of \cite{Get94} by our present consideration. Namely, one can first verify the derivation property. Then all the remaining properties can be checked on the cup product generators $N^0$ and $1 \otimes \frm \subset N^1$. In particular, the fact that $d'$ and $d''$ lift the corresponding ones under \eqref{OCar_hG_to_hg} will be obvious. As a consequence, to verify $[d,\delta]=0$ on generators, if the target lands on $N^0$, then it automatically holds, since the map \eqref{OCar_hG_to_hg} is an isomorphism on $N^0$. 
From this, we see that the only non-trivial things to check is that
\eqref{dd_eqn2} is consistent with it being a derivation, and that $[d'',\delta''] + [d',\delta'] = 0$ on $1 \otimes \frm \subset N^1$. For both of these verifications, we have found it more convenient to use a more classical notation close to the one in \cite{Get94}. Namely, for a $\cO(G)$-comodule $V$, we think of $C^n(G,V)$ as the set of algebraic maps $G^n \ra V$. For $F \in C^n(G,V)$ and for $1 \leq i \leq n$, define $d_e^{(i)}(F) \in C^{n-1}(G,V \otimes \widetilde{\frg^*})$ by
\begin{equation*}
d_e^{(i)}(F)(g_1,\ldots,g_{n-1}) \, = \, 
(\id_V \otimes \Ad^*_{g_1\ldots,g_{i-1}})(d_g(F(g_1,\ldots,g_{i-1},g,g_{i},\ldots,g_{n-1}))|_{g=e}))
\end{equation*}
where $\Ad^*$ is the (left) coadjoint action of $G$ on $\frg^*$.
Then define $d'' : N^n(G,\OCar^*(\frg^*,A)) \ra N^{n-1}(G,\OCar^*(\frg^*,A))$ by 
$d''(F) = \sum_{i=1}^{n} (-1)^{i-1} \overline{d_e^{(i)}(F)}$, where the overline refers to the map $\OCar^p(\frg^*,A) \otimes \widetilde{\frg^*} \ra \OCar^{p+1}(\frg^*,A)$ given by multiplication.
For $F_1 \in N^p(G,\OCar^*(\frg^*,A))$ and $F_2 \in C^q(G,\OCar^*(\frg^*,A))$, we can then plug in \eqref{cup_OYG_2} to see that
$\overline{d^{(i)}_e(F_1\cup F_2)} = \overline{d^{(i)}_e(F_1)}\cup F_2$ for $1 \leq i \leq p$ and $\overline{d^{(p+j)}_e(F_1\cup F_2)} = F_1\cup \overline{d^{(j)}_e( F_2)}$
for $1 \leq j \leq q$ (this first equation requires that $F_1 \in N^p$ but the second equation is also true for $F_1 \in C^p$), which implies the derivation property.
The verification that $[d'',\delta''] = [d',\delta'] = 0$ on $1 \otimes \frm \subset N^1$ is also straightforward from this formula of $d''$.
%
%
%
\epf

Thus, if we collapse the first and second equivariant degrees (for each fixed Hodge degree there are only finitely many first equivariant degrees, so there is no ambiguity in this collapsing), and take $\delta = \delta' + \delta''$ and $d = d'+d''$, then $N^{\bullet}(G,\OCar^*(\frg^*,A))$ becomes a tri-complex.

\bdf
The tri-complex $(N^{\bullet}(G,\OCar^*(\frg^*,A)),\partial,\delta,d)$ is called the \emph{Cartan-Getzler-de Rham tri-complex}.
\edf

One can then follow the same procedure to define shifted symplectic structures based on the Cartan-Getzler-de Rham tri-complex. In particular, we may form the bicomplex $\DR^*_{\CG}(Y/G)$ by collapsing the equivariant degree into the intrinsic one. Here, we take the direct product total complex
\begin{equation*}
\DR^p_{\CG}(Y/G) \, := \, \OCar^p(\frg^*,A)^{\hG} \, := \, \Tot^{\Pi}\, N^{\bullet}(G,\OCar^p(\frg^*,A))
\end{equation*}
and likewise define $\DRcl_{\CG}(Y/G)$ and $F^n\DRcl_{\CG}(Y/G)$ as above. Notice that, in this case, a cocyle $\omega \in \DR^2_{\CG}(Y/G)$ may not define a $G$-equivariant map \eqref{omega_sharp}. However, there is always a map of complexes $\DR^2_{\CG}(Y/G) \ra \OCar^2(Y/G)$, so that $\omega$ still defines a map \eqref{omega_sharp} of DG-modules over $A$. 

\bdf  \label{shifted_sympl_CG_def}
Assume that $A$ is almost cofibrant.
Then an \emph{$m$-shifted symplectic structure in Cartan-Getzler model} on $[Y/G]$ is a cocycle $\widetilde{\omega} \in Z^{m+2} F^2\DRcl_{\CG}(Y/G)$ whose induced map \eqref{omega_sharp} is a quasi-isomorphism.
One can define spaces of (pre-)symplectic structures in Cartan-Getzler model in the same way as above.
\edf

\brm
One may try to construct an analogue of Getzler's extension for polyvector fields. In fact, there is always a map of bicomplexes (see \eqref{Hom_hG_to_CE})
\begin{equation}  
N(G,\XCar^p(\frg,A,m)) \raq
\CE(\frg^*, \XCar^p(\frg,A,m)) 
\end{equation}
and we might again try to lift the brackets $\{-,-\}'$ and $\{-,-\}''$ to $N(G,\XCar^p(\frg,A,m))$. However, this seems to be more complicated than the case for differential forms. For example, since the cup product is noncommutative, it seems difficult to write a formula for $\{-,-\}'$ that is anti-symmetric.
\erm

\section{Functoriality of Cartan polyvectors}  \label{sec_functoriality}

In this section, we discuss the functoriality of our previous constructions. A map $[Y'/G'] \ra [Y/G]$ of global quotients will mean a pair of maps $Y' \ra Y$ and $G' \ra G$ commuting with the actions.
It is clear that all of our constructions involving Cartan forms (as well as Getzler's extension) are functorial with respect to such maps. 
For Cartan polyvectors, there is also a functoriality in terms of a span.
We first consider the non-equivariant case. For any map $f:A \ra B$ of CDGA's, there are the obvious maps
\begin{equation}  \label{Pol_cospan_1}
\Pol^p(A,m) \xraq{f_*}
\Homcom_{A}(\, \Sym^p_{A}(\Omega^1(A)[m+1]) \, , \, B \,) 
\xlaq{f^*} 
\Pol^p(B,m)
\end{equation}
which allows one to define
\begin{equation} \label{Pol_rel_1}
\Pol^p(f,m) := \{ \, (\theta_A,\theta_B) \in \Pol^p(A,m) \times \Pol^p(B,m) \, | \, f_*(\theta_A)=f^*(\theta_B)  \, \}
\end{equation}
which is an $(m+1)$-shifted DG Lie subalgebra of $\Pol^p(A,m) \times \Pol^p(B,m)$. Indeed, the corresponding statement holds at the complexes $C^p(A,m) := \Homcom_k((A[m+1])^{\otimes p},A)$. Namely, for any many $f : A \ra B$ of graded vector spaces, one can take
\begin{equation*}
C^p(A,m)  \xraq{f_*}  \Homcom_k((A[m+1])^{\otimes p},B)  \xlaq{f^*}  C^p(B,m)
\end{equation*}
and use it to form $C^p(f,m)$ in the same way. The subcomplex $C^p(f,m)[m+1] \subset C^p(A,m)[m+1] \times C^p(B,m)[m+1]$ is then preserved by the pre-Lie structure \eqref{Hoch_pre_Lie_1}. This implies that the Schouten-Nijenhuis bracket preserves the subcomplex $\Pol^p(f,m) \subset \Pol^p(A,m) \times \Pol^p(B,m)$.

Thus, we have maps of DG Lie algebras
\begin{equation*}
\Pol^p(A,m)[m+1] \xlaq{\pi_A} \Pol^p(f,m)[m+1] \xraq{\pi_B} \Pol^p(B,m)[m+1] 
\end{equation*}

As a cochain complex, $\Pol^p(f,m)$ is defined as a certain kernel. To ensure that it has good homological property, we would like to ensure that this kernel is quasi-isomorphism to the fiber (or cocone). We will use a sufficient condition that is different from the one used in \cite{Pri17}:

\blm  \label{retract_then_surj}
Suppose that the map $\Omega^1(A)\otimes_{A} B \ra \Omega^1(B)$ has a retract (\ie, left inverse) as a map of graded $B$-modules, then the map 
$f^* : \Pol^p(B,m) \ra \Homcom_A( \Sym^p_{A}(\Omega^1(A)[m+1])  ,  B ) $ is surjective.
\elm

Now we discuss the equivariant case. Thus, consider a map $\pi:[Y'/G'] \ra [Y/G]$ of global quotients for $Y = \Spec\, A$ and $Y' = \Spec\, A'$. Then we have the following diagram:
\begin{equation}  \label{quot_map_cospan}
\begin{tikzcd}[column sep = 1.5em]
\XCar^p(\frg,A,m)^G \ar[r, "(1)"] \ar[d, hook]
& \Homcom_{A}( \Sym^p_{A}(\OCar^1(\frg^*,A)[m+1])  ,  A' )^{G'}  \ar[d, hook]
& \XCar^p(\frg',A',m)^{G'} \ar[l, "(2)"'] \ar[d, hook] \\
\CE(\frg^{*},\XCar^p(\frg,A,m)) \ar[r] \ar[d, equal]
& \CE(\frg^{\prime *},\Homcom_{A}( \Sym^p_{A}(\OCar^1(\frg^*,A)[m+1]) , A')) \ar[d, equal]
& \CE(\frg^{\prime *}, \XCar^p(\frg',A',m)) \ar[l] \ar[d, equal] \\
\Pol^p(B,m) \ar[r]
& \Homcom_{B}( \Sym^p_{B}(\Omega^1(B)[m+1]) , B') 
& \Pol^p(B',m) \ar[l, "(3)"']
\end{tikzcd}
\end{equation}
where the last line is the corresponding maps \eqref{Pol_cospan_1} for the map $\pi^{\sharp}_{\CE} : B := \cO([Y/\frg]) \ra \cO([Y'/\frg']) =: B'$. In the identification between the second and third row, we implicitly assert (with a completely parallel proof) an analogue of Proposition \ref{Pol_CE_isom} for this relative case.

Collapse the equivariant and intrinsic degrees, and form the analogue of \eqref{Pol_rel_1} for each row, then we have maps of $(m+1)$-shifted weight graded DG Lie algebras
\begin{equation}  \label{quot_map_span}
\begin{tikzcd}
\PCar^*(Y/G,m)  \ar[d, hook]
& \PCar^*(\pi,m)  \ar[d, hook] \ar[l, "(4)"'] \ar[r, "(5)"]
& \PCar^*(Y'/G',m) \ar[d, hook] \\
\Pol^*(B,m) 
& \Pol^*(\pi^{\sharp}_{\CE},m) \ar[l, "(6)"'] \ar[r]
& \Pol^*(B',m) 
\end{tikzcd}
\end{equation}

We will now focus on the special case of a triangular extension, in the sense of the following

\bdf  \label{triang_ext_def}
A map $\pi:[Y'/G'] \ra [Y/G]$ of global quotients is said to be a \emph{triangular extension} if there exists isomorphisms $Y' \cong Y \times H$ and $G' \cong G \times H$ for some linear algebraic group $H$, such that under this identification, we have
\begin{enumerate}
	\item $H$ acts on $H$ by right translation.
	\item The map $\pi:[Y'/G'] \ra [Y/G]$ is the obvious projections to $Y$ and $G$.
\end{enumerate}
Notice that (2) implies in particular that $H$ acts trivially on $Y$, and $G$ acts on $Y$ the same way as in $[Y/G]$. 
\edf

\blm  \label{triang_ext_lem}
Let $\pi : [Y'/G'] \ra [Y/G]$ be a triangular extension. Then the map $\OCar^1(Y/G) \otimes_{A} A' \ra \OCar^1(Y'/G')$ is a quasi-isomorphism.
\elm

\bpf
The map is injective, whose quotient is the cone of the identity map on the free module $\frh^* \otimes A'$.
\epf

\blm  \label{prin_bundle_equiv}
Let $\pi : [Y'/G'] \ra [Y/G]$ be a triangular extension. For any $M,N \in \dgMod_{G}(Y)$, consider the functoriality maps 
\begin{equation*}
\begin{split}
\Homcom_A(M,N)^G &\raq \Homcom_{A'}(\pi^*M, \pi^*N)^{G'}
\\
\Tot^{\Pi}_{\Delta} \, C^{\bullet}(G,\Homcom_A(M,N)) &\raq \Tot^{\Pi}_{\Delta} \, C^{\bullet}(G',\Homcom_{A'}(\pi^*M, \pi^*N))
\end{split}
\end{equation*}
Then the first map is an isomorphism; while the second is a quasi-isomorphism.
\elm

\bpf
Notice that $\pi: Y' \ra Y$ is a (trivial) principal $H$-bundle, hence the functor $\pi^*$ and $\pi_*(-)^H$ are inverse equivalence between the categories $\dgMod(\cO(Y))$ and $\dgMod_H(\cO(Y'))$. One can remember the $G$-equivariance structures when applying these two functors to obtain an equivalence of categories between $\dgMod_{G}(\cO(Y))$ and $\dgMod_{G'}(\cO(Y'))$. The first statement then follows from the obvious DG enhancements of both of these functors. 

For the second statement, notice that the cosimplicial complex in the right hand side is the diagonal of a bi-cosimplicial complex
$C^{\bullet}(G, C^{\bullet}(H,\Homcom_{A'}(\pi^*M, \pi^*N)))$. 
We may therefore first take the homotopy limit with respect to $C^{\bullet}(H,-)$. Since $\pi^*N$ is a trivial extension with respect to $H$, it has split diagonal as an $H$-module, so that by Corollary \ref{split_diag_equiv}, we have a quasi-isomorphism $\Homcom_{A'}(\pi^*M, \pi^*N)^H \xra{\simeq} C^{\bullet}(H,\Homcom_{A'}(\pi^*M, \pi^*N))$. By the same argument as the above paragraph, we have an isomorphism $\Homcom_{A'}(\pi^*M, \pi^*N)^H \cong \Homcom_A(M,N)$.
In other words, we have a quasi-isomorphism $\Homcom_A(M,N) \xra{\simeq} \Tot^{\Pi}_{\Delta} \, C^{\bullet}(H,\Homcom_{A'}(\pi^*M, \pi^*N))$, which also extends degreewise to $C^{\bullet}(G,-)$. Taking the homotopy limit, we obtain the desired quasi-isomorphism.
\epf

\bpp  \label{triang_ext_functoriality}
Let $\pi : [Y'/G'] \ra [Y/G]$ be a triangular extension. Assume that $G'$ is reductive, and $\Omega^1(Y)$ is cofibrant. Then in \eqref{quot_map_cospan}, the map (1) is an isomorphism, while the maps (2),(3) are surjective quasi-isomorphisms after collapsing the equivariant and intrinsic gradings.
As a consequence, in \eqref{quot_map_span}, the map (5) is an isomorphism, while the maps (4), (6) are quasi-isomorphisms.
\epp

\bpf
The fact that (1) is an isomorphism follows from Lemma \ref{prin_bundle_equiv}.
To show that (2) and (3) are quasi-isomorphisms, notice that under our present assumption both of the $A'$-modules in the quasi-isomorphism of Lemma \ref{triang_ext_lem} are cofibrant. Taking the $A'$-dual of their shifted symmetric powers therefore remains a quasi-isomorphism. Thus, we conclude that (2) is a quasi-isomorphism by reductiveness of $G'$ (note that the exactness of taking $G$-invariants holds not just for rational modules, but also for Hom complexes between them, see, \eg, Proposition \ref{Reynold_prop}).
Similarly, to conclude that (3) is a quasi-isomorphism, simply consider the second row. 
To see that (2) and (3) are surjective, one may use the same observation as in Lemma \ref{retract_then_surj}, applied to $\OCar^1$.
The statements for the maps (4),(5),(6) in \eqref{quot_map_span} are then formal consequences.
\epf

More generally, we may consider any sequence of maps of global quotients
\begin{equation*}
\pi \, = \, \left[ \, [Y_0/G_0] \xla{\pi_1} [Y_1/G_1] \xla{\pi_2} \ldots \xla{\pi_r} [Y_r/G_r] \, \right]
\end{equation*}
and we may form 
\begin{equation*}  
\begin{split}
\PCar^p(\pi,m) \, &:= \, {\rm eq} \, \left[ \, \Pi_{i=0}^r \PCar^p(Y_i/G_i,m) \, \,
\substack{\rightarrow \\ \rightarrow} \, \,
\Pi_{i=0}^{r-1} \, \Homcom_{A_i}( \Sym_{A_i}^p(\OCar^1(Y_i/G_i)[m+1]), A_{i+1})^{G_{i+1}} \, \right] \\
\Pol^p(\pi^{\sharp}_{\CE},m) \, &:= \, {\rm eq} \, \left[ \, \Pi_{i=0}^r \Pol^p(B_i,m) \, \,
\substack{\rightarrow \\ \rightarrow} \, \,
\Pi_{i=0}^{r-1} \, \Homcom_{B_i}( \Sym_{B_i}^p(\Omega^1(B_i)[m+1]), B_{i+1}) \, \right]
\end{split}
\end{equation*}
where we write $B_i := \cO([Y_i/\frg_i])$ in the second line.
Both of these are then (shifted) DG Lie subalgebras of the product $\Pi_{i=0}^r \PCar^*(Y_i/G_i,m)$ and $\Pi_{i=0}^r \Pol^p(B_i,m)$ respectively. A repeated application of \eqref{triang_ext_functoriality} then gives the following analogue of \cite[Properties 2.5]{Pri17}:
\bpp  \label{length_r_Pol}
Assume that each $\pi_i$ is a triangular extension. Then we have
\begin{enumerate}
	\item There is a map of DG Lie algebras $\PCar^*(\pi,m)[m+1] \ra \Pol^p(\pi^{\sharp}_{\CE},m)[m+1]$ commuting with the projections to each component.
	\item Denote by $\cP(\pi)$ either $\PCar^*(\pi,m)[m+1]$ or $\Pol^p(\pi^{\sharp}_{\CE},m)[m+1]$. Then the map
	\begin{equation*}
	\cP(\pi) \raq \cP(\pi_1,m) \times^h_{\cP([Y_1/G_1])} \ldots \times^h_{\cP([Y_{r-1}/G_{r-1}])} \cP(\pi_r,m)
	\end{equation*} 
	is a quasi-isomorphism of DG Lie algebras.
	\item The projection $\PCar^*(\pi,m) \ra \PCar^*([Y_i/G_i],m)$ is a quasi-isomorphism for each $0 \leq i \leq r$.
	\item The projection
	$\Pol^*(\pi^{\sharp}_{\CE},m) \ra \PCar^*([Y_0/\frg_0],m)$ is a quasi-isomorphism.
\end{enumerate}
\epp

\section{Comparison with the general notions} \label{comparison_sec}

In this section, we show that for global quotients, our notions of shifted symplectic and Poisson structures coincides with that in \cite{Pri17} (see Theorems \ref{shifted_sympl_equiv} and \ref{shifted_Poiss_equiv} below).

Given any simplicial derived affine scheme $X : \Delta^{\op} \ra \dAff$, we may thicken it into a bisimplicial derived affine scheme $\bX$ by taking
\begin{equation}  \label{bXij}
\bX_{i,j} \, := \, \Hom_{\setdel}(\Delta^i \times \Delta^j , X)
\end{equation}
understood as an end. In the notation of \cite{Pri17}, this is $\bX_{i,\bullet} = X^{\Delta^i}$. We will call $i$ the horizontal degree and $j$ the vertical degree. Of course, the horizontal and vertical degrees are symmetric for now, but we will apply asymmetric operations on them later on.

We will focus on the case $X = [Y/G]_{\bullet}$. As a $1$-groupoid, any map from a simplicial set to $X$ is determined by its $1$-skeletal part, 
which is required to satisfy a functoriality relation specified at its $2$-skeletal part. Thus, an element in $\bX_{i,j}$ is the same as a commutative $(i \times j)$-square in the groupoid $G^{\op} \ltimes Y$ (see our convention \eqref{Y_mod_G_def}). From this it is easy to see (\cf \, \cite[Example 3.6]{Pri17}) that $\bX_{m,\bullet}$ is isomorphic to $[Y \! \times \! G^m/G^{m+1}]_{\bullet}$ with action given by
\begin{equation*}
 (y,h_1,\ldots,h_m) \cdot (g_0,\ldots,g_m) = 
(y g_0, g_0^{-1} h_1 g_1, \ldots,  g_{m-1}^{-1} h_m g_m)
\end{equation*} 
It is clear from this description that, for any $[n]\ra [m]$ in $\Delta$, the map $\bX_{m,\bullet} \ra \bX_{n,\bullet}$ of simplicial derived affine schemes is induced by a map of global quotients $[Y \! \times \! G^m / G^{m+1}] \ra [Y \! \times \! G^n / G^{n+1}]$. Namely the one $Y \times G^m \ra Y \times G^n$ from the simplicial structure of $X_{\bullet}$, and the obvious map $G^{m+1} \ra G^{n+1}$ induced from the powering indices.
We will write
\begin{equation*}
[X_n / G_n] \, := \, [Y \! \times \! G^n / G^{n+1}]
\end{equation*}
The following Lemma is clear:
\blm  \label{inj_triang_ext}
For any monomorphism $[n] \rinto [m]$ in $\Delta$, the map of global quotients $[X_m/G_m] \ra [X_n / G_n]$ is a triangular extension in the sense of Definition \ref{triang_ext_def}.
\elm

Recall the functor $\Gamma^* : \CAlg_k^{\Delta} \ra \CDGAp_k$ from Appendix \ref{app_small_alg} (see the discussion preceding Proposition \ref{Gamma_star_1thin}), or rather its bigraded version.
Apply it to each $\cO(\bX_{m,\bullet})$, and denote the result simply as $\Gamma^* \cO(\bX)$, with the understanding that $\Gamma$ and $\Gamma^*$ always concerns the vertical degree (\ie, applied to each fixed horizontal degree). Denote by $X^{\vc}$ the ``vertically constant'' bisimplicial derived affine scheme $X^{\vc}_{i,j} := X_i$, then the map $\Delta^j \ra \Delta^0$ induces by \eqref{bXij} a map $X^{\vc} \ra \bX$. Applying $\Gamma \Gamma^*$ to $\cO(\bX) \ra \cO(X^{\vc})$, and we obtain
\begin{equation*}
\cO(\bX) \raq \Gamma \Gamma^* \cO(\bX) \raq \cO(X^{\vc})
\end{equation*}
We recall the following result from \cite{Pri17}:

\bpp  \label{bX_GGX_X_hoeq}
Both of the maps $X^{\vc} \raq \Spec \, \Gamma \Gamma^* \cO(\bX) \raq \bX$ of bisimplicial derived affine schemes are simplicial homotopy equivalences at each fixed vertical degree. 
\epp

\bpf
The statement for the first map is the content of (the proof of) \cite[Proposition 3.13]{Pri17}. At a fixed vertical degree, the composition is the canonical map $X \ra X^{\Delta^j}$, which is also a simplicial homotopy equivalence since the map $\Delta^j \ra \Delta^0$ of simplicial sets is.
\epf

As noted by \cite[Example 3.6]{Pri17} (see also Proposition \ref{Gamma_star_1thin} and \ref{glob_quot_1thin_quot}), 
we have $\Gamma^*(\cO([Y/G]_{\bullet})) \cong \cO([Y/\frg])$.
Thus, in the case $X = [Y/G]_{\bullet}$, the bigraded CDGA $\Gamma^* \cO(\bX_{m,\bullet})$ is given by
\begin{equation*}
\Gamma^* \cO(\bX_{m,\bullet}) \, \cong \, \cO([X_m/\frg_m])
\end{equation*}

In \cite{Pri17}, shifted symplectic and Poisson structures on $X$ are defined in terms of the differential calculus on the cosimplicial bigraded CDGA $\Gamma^* \cO(\bX)$. In the case $X = [Y/G]_{\bullet}$, we may rewrite these in terms of the Cartan forms and Cartan polyvectors, in view of Propositions \ref{CE_Omega_Car_Omega_CE} and \ref{Pol_CE_isom}. We will now assume $X = [Y/G]_{\bullet}$ throughout.

Let $\OCar^p(\bX)$ be the DG module over the bi-cosimplicial CDGA $\cO(\bX)$ such that for each fixed horizontal degree $m$, $\OCar^p(\bX)_{m,\bullet}$ is the strictly Cartesian module corresponding to $\OCar^p(X_m/G_m)$. In other words,
\begin{equation}  \label{OCar_bX}
\OCar^p(\bX)_{m,n} \, = \, C^n(G_m,\OCar^p(X_m / G_m))
\end{equation}
then we have
\blm  \label{OCar_bX_hcart}
The DG module $\OCar^p(\bX)$  is strictly Cartesian in the vertical direction, and homotopy Cartesian in the horizontal direction.
\elm

\bpf
Strict Cartesianness in the vertical direction holds by construction. To verify homotopy Cartesianness in the horizontal direction, recall by Remark \ref{check_at_di} that it suffices to check at injective maps $[m] \rinto [m']$. Since the corresponding maps of global quotients is a triangular extension (see Lemma \ref{inj_triang_ext}), the result holds by Lemma \ref{triang_ext_lem}. 
\epf

In view of Propositions \ref{Gamma_star_1thin} and \ref{glob_quot_1thin_quot}, taking the Chevalley-Eilenberg complex amounts to the passage to the $1$-thin quotient $\Gamma\Gamma^* \cO(\bX)$. Thus, Propositions \ref{CE_Omega_Car_Omega_CE} gives an isomorphism
\begin{equation}  \label{OCar_GGstar}
\Tot \, N( \OCar^p(\bX) \otimes_{\cO(\bX)} \Gamma \Gamma^*( \cO(\bX) ) ) \, \cong \, \Tot \, \Omega^p(\Gamma^* \cO(\bX))
\end{equation}
where we apply the normalization $N$ to the vertical direction (\ie, the said isomorphism holds for each fixed horizontal degree).
Notice that on both sides, since we have fixed the Hodge degree $p$, the bicomplex is bounded in the equivariant degree, and so there is no ambiguity in taking the total complex.

Similarly, in view of the discussion at the paragraph of \eqref{Hom_hG_to_CE_tensor}, Proposition \ref{Pol_CE_isom} gives an isomorphism
\begin{equation}  \label{XCar_GGstar}
\Tot \, N( \Homcom_{\cO(\bX)}^{\sim}( \, \Sym^p_{\cO(\bX)}( \OCar^1(\bX)[m+1] ) \, , \, \Gamma \Gamma^*(\cO(\bX)) \,) ) \, \cong \, \Tot \, \Pol^p(\Gamma^* \cO(\bX),m)
\end{equation}
also understood to be taken  for each fixed horizontal degree, and there is also no ambiguity in taking the total complex.

\brm
Notice that the isomorphisms \eqref{OCar_GGstar} and \eqref{XCar_GGstar} do not hold before taking the total complex. Indeed, in the definition \eqref{OCar_bX} of $\OCar(\bX)$ (and hence the left hand sides of \eqref{OCar_GGstar} and \eqref{XCar_GGstar}), we have already collapsed the first equivariant grading with the intrinsic grading. On the other hand, on the right hand side of \eqref{OCar_GGstar} and \eqref{XCar_GGstar}, the first equivariant grading is collapsed with the second equivariant grading, so that we have different double complexes.
\erm

For both \eqref{OCar_GGstar} and \eqref{XCar_GGstar}, there are canonical maps from the versions before passing to the $1$-thin quotient $\Gamma\Gamma^* \cO(\bX)$. Thus, there are maps 
\begin{equation}  \label{OCar_bX_three_terms}
\OCar^p(X_n/G_n)^{G_n} \raq \Tot^{\Pi} \, N(\OCar^p(\bX_{n,\bullet})) \raq \Tot \, \Omega^p(\Gamma^* \cO(\bX_{n,\bullet}))
\end{equation}
\begin{equation} \label{XCar_bX_three_terms}
\XCar^p(X_n/G_n,m)^{G_n} \raq 
\Tot^{\Pi} \, N(G, \XCar^p(X_n/G_n,m)) \raq \Tot \, \Pol^p(\Gamma^* \cO(\bX_{n,\bullet}),m)
\end{equation}
where we recall that the $\Tot$ in the third term of both \eqref{OCar_bX_three_terms} and \eqref{XCar_bX_three_terms} are equal to their $\Tot^{\Pi}$, so that there is indeed a map from $\Tot^{\Pi}$ of the middle terms. We also recall that the second map in \eqref{XCar_bX_three_terms}
is given by the normalization of \eqref{Hom_hG_to_CE_tensor}.

Notice that each term of \eqref{OCar_bX_three_terms} is part of a de Rham complex. Namely, if we consider the global quotient $[X_n/G_n]$, then the three terms in \eqref{OCar_bX_three_terms} are respectively the Hodge degree $p$ part of the Cartan-de Rham bicomplex, the Cartan-Getzler-de Rham bicomplex, and the ordinary de Rham complex of $\cO([X_n / \frg_n])$. The maps in \eqref{OCar_bX_three_terms} also commutes with the respective de Rham differentials.
\bpp  \label{OCar_bX_GGOX_holim}
Assume that $A$ is almost cofibrant, then the map
\begin{equation*}
\holim_{[n] \in \Delta} \, \Tot^{\Pi} \, N(\OCar^p(\bX_{n,\bullet})) \raq \holim_{[n] \in \Delta} \, \Tot \, \Omega^p(\Gamma^* \cO(\bX_{n,\bullet}))
\end{equation*}
induced by \eqref{OCar_bX_three_terms} is a quasi-isomorphism.
\epp

\bpf
Recall from Proposition \ref{bX_GGX_X_hoeq} that the map $i : \Spec \, \Gamma \Gamma^*(\cO(\bX)) \ra \bX$ is a simplicial homotopy equivalence for each fixed vertical degree, so that in particular, they have equivalent $\infty$-categories of quasi-coherent sheaves (defined in terms of homotopy Cartesian modules, see Appendix \ref{app_cart}). By Lemma \ref{OCar_bX_hcart}, the object $\OCar^p(\bX)$ lies in $\dgMod(\cO(\bX))_{\hca}$. The left hand side is the derived Hom complex between $\cO(\bX)$ and $\OCar^p(\bX)$ inside $\dgMod(\cO(\bX))_{\hca}$, while the right hand side is the derived Hom complex between $i^*\cO(\bX)$ and $i^*\OCar^p(\bX)$ inside $\dgMod(\Gamma \Gamma^* \cO(\bX))_{\hca}$.
\epf

As a consequence, we have the following
\bthm  \label{shifted_sympl_equiv}
Assume that $A$ is almost cofibrant, then the space of shifted symplectic structures in Cartan-Getzler model is homotopy equivalent to the space of shifted symplectic structures as defined in \cite{Pri17}.
If $G$ is moreover reductive, then it is also homotopy equivalent to the space of shifted symplectic structures in Cartan model.
\ethm

\bpf
Clearly, the second statement follows from the first. 

Notice that the left hand side $\Tot^{\Pi}\,N(\OCar^p(\bX_{n,\bullet}))$ of Proposition \ref{OCar_bX_GGOX_holim} is constant as a simplicial system in $n$. In other words, for each $[n] \ra [n']$, the induced map is a quasi-isomorphism. 
Indeed, it suffices to show this for a coface map, for which the maps of total quotient is a triangular extension (see Lemma \ref{inj_triang_ext}), so that the statement follows from Lemmas \ref{triang_ext_lem} and \ref{prin_bundle_equiv}.
Thus, the projection from homotopy limit on the left hand side of Proposition \ref{OCar_bX_GGOX_holim} to the term at $n = 0$ is a quasi-isomorphism.
As a result, Proposition \ref{OCar_bX_GGOX_holim} implies that the spaces of pre-symplectic structures in Cartan-Getzler model and in \cite{Pri17} are homotopy equivalent, so it suffices to verify that the notion of non-degeneracy matches. 
By Cartesianness, non-degeneracy%
\footnote{In \cite{Pri17}, one first defines the spaces ${\rm PreSp} \subset {\rm Sp}$ at each cosimplicial degree, and then takes the homotopy limit. Since ${\rm PreSp} \subset {\rm Sp}$ is a union of connected components, the same is true for their homotopy limits, \ie, non-degeneracy is a property of a pre-symplectic structure, which should be checked at the projection to each $[n] \in \Delta$, but by Cartesianness suffices to be checkd at $[0]$.} 
of $\omega \in \holim_{[n] \in \Delta} \, \Tot \, \Omega^2(\Gamma^* \cO(\bX_{n,\bullet}))$ also depends only on its projection to the $n=0$ term. By \cite[Definition 3.21]{Pri17} it is defined by first taking $-\otimes_{\cO[Y/\frg]} \cO(Y)$, which precisely coincides with our notion of non-degeneracy in the Cartan-Getzler model.
\epf

Next, we compare the notions of shifted Poisson structures.
Note that the first and third terms of \eqref{XCar_bX_three_terms} are the weight $p$ parts of the $m$-shifted Cartan polyvectors of $[X_n/G_n]$ and the $m$-shifted polyvectors of $\cO([X_n/\frg_n])$ respectively. Hence, this is a map of $(m+1)$-shifted DG Lie algebras.
In some sense, we would like to relate the two by taking the homotopy limit of \eqref{XCar_bX_three_terms} over $[n] \in \Delta$. Of course, the terms in \eqref{XCar_bX_three_terms} is not really functorial in $[n] \in \Delta$, but in Section \ref{sec_functoriality}, we have developed enough of an $\infty$-categorical functoriality for our purposes. 
Namely, we first recall the following result (see, \eg, \cite[Example 8.5.12]{Rie} for a proof).

\blm  \label{Delta_inj_ho_init}
Let $\Delta_{{\rm inj}} \subset \Delta$ be the subcategory consisting of all objects, and only the monomorphisms. Then the inclusion $\Delta_{{\rm inj}} \rinto \Delta$ is homotopy initial.
\elm

From now on, assume that $A$ is almost cofibrant and $G$ is reductive.
Above a map $[n] \ra [n']$ in $\Delta_{{\rm inj}}$, the global quotient $[X_{n'}/G_{n'}] \ra [X_n/G_n]$ is a triangular extension (see Lemma \ref{inj_triang_ext}), and hence, by Proposition \ref{triang_ext_functoriality}, gives a weak functoriality in terms of a span
\eqref{quot_map_span} of DG Lie algebras. The coherent functoriality of these spans is established in Proposition \ref{length_r_Pol}, so that we may regard both sides of \eqref{XCar_bX_three_terms} the terms at $[n]$ of $\infty$-functors $\Delta_{{\rm inj}} \ra {\it DGLA}^{gr}_k$ to the $\infty$-category of weight graded DG Lie algebras (where the bracket has weight $-1$). 
Taking the $\infty$-categorical limits, we have a map of weight graded DG Lie algebras
\begin{equation}   \label{infty_lim_Pol}
\underset{[n] \in \Delta_{{\rm inj}}}{{\rm lim}^{\infty}} \, \Pol_{{\rm Car}}^*(X_n / G_n,m)[m+1] \raq 
\underset{[n] \in \Delta_{{\rm inj}}}{{\rm lim}^{\infty}} \, \Tot \, \Pol^*(\cO([X_n / \frg_n],m)[m+1]
\end{equation}

\blm  \label{infty_lim_Pol_RHom}
If we forget about the Lie structures, and only consider the underlying cochain complexes, then the two sides of \eqref{infty_lim_Pol} are given respectively by
\begin{equation*}
\begin{split}
\underset{[n] \in \Delta_{{\rm inj}}}{{\rm lim}^{\infty}} \, \Pol_{{\rm Car}}^p(X_n / G_n,m) 
\, &\simeq \, 
{\bm R}\Homcom_{\cO(\bX)|_{\Delta_{{\rm inj}}\times \Delta}}( \Sym_{\cO(\bX)}^p(\OCar^1(\bX) [m+1]),\cO(\bX))
\\
\underset{[n] \in \Delta_{{\rm inj}}}{{\rm lim}^{\infty}} \, \Tot \, \Pol^p(\cO([X_n / \frg_n]),m) 
\, &\simeq \, 
{\bm R}\Homcom_{\Gamma \Gamma^* \cO(\bX)|_{\Delta_{{\rm inj}}\times \Delta}}( i^* \Sym_{\cO(\bX)}^p(\OCar^1(\bX) [m+1]), i^*\cO(\bX))
\end{split}
\end{equation*}
where $i^* : \dgMod(\cO(\bX)) \ra \dgMod(\Gamma \Gamma^* \cO(\bX))$ is the tensor map.
\elm

\bpf
To compute the right hand side, we may apply \eqref{RHom_inflim_Homsim} to $I = \Delta_{{\rm inj}}$, so that the computation is reduced to individual ${\bm R}\Hom_{\cO(\bX_{n,\bullet})}(-,-)$, which are precisely given by the left hand side by applying \eqref{Homhc_holim_Homsim}, since the modules in question are strictly Cartesian in the vertical direction.
\epf

As a result, we have
\bthm  \label{shifted_Poiss_equiv}
Assume that $A$ is almost cofibrant and $G$ is reductive, then the map \eqref{infty_lim_Pol} of DG Lie algebras is a quasi-isomorphism. Hence, the space of $m$-shifted Poisson structures in Cartan model is homotopy equivalent to the space of $m$-shifted Poisson structures as defined in \cite{Pri17}.
\ethm

\bpf
By Lemma \ref{infty_lim_Pol_RHom}, we have written both sides of \eqref{infty_lim_Pol} as derived Hom complexes between homotopy Cartesian DG modules. By Lemma \ref{Delta_inj_ho_init}, we have equivalences of $\infty$-categories ${\it dgMod}(\cO(\bX))_{\hca} \simeq {\it dgMod}(\cO(\bX)|_{\Delta_{{\rm inj}}\times \Delta})_{\hca}$ and ${\it dgMod}(\Gamma \Gamma^* \cO(\bX))_{\hca} \simeq {\it dgMod}(\Gamma \Gamma^* \cO(\bX)|_{\Delta_{{\rm inj}}\times \Delta})_{\hca}$ since they can be written as a homotopy limit (see Appendix \ref{app_cart}). By Proposition \ref{bX_GGX_X_hoeq}, we also have an equivalence ${\it dgMod}(\cO(\bX))_{\hca} \simeq {\it dgMod}(\Gamma \Gamma^* \cO(\bX))_{\hca}$. In particular, they have quasi-isomorphic ${\bm R}\Hom$.

To prove the second statement, notice that, by Proposition \ref{triang_ext_functoriality}, the $\infty$-functor on the left hand side of \eqref{infty_lim_Pol} is in fact a constant functor, so that the projection to the term at $[0] \in \Delta_{{\rm inj}}$ is a quasi-isomorphism. This term is precisely our Cartan model for polyvectors.
\epf

\appendix

\section{Cartesian modules and quasi-coherent sheaves}  \label{app_cart}

\bdf
Given a small category $I$, and a functor $A : I \ra \DGA_k$, then a DG module over $A$ (often simply called an $A$-module) consists of a (right) DG module $M_i$ over $A_i$ for each $i \in \Ob(I)$, together with a map $\varphi : M_i \ra M_j$ of DG modules, linear over $\varphi : A_i \ra A_j$, for each $\varphi \in \Hom_I(i,j)$, satisfying the obvious functoriality. Denote by $\dgMod(A)$ the category of DG modules over $A$. This is the underlying category of a DG category $\dgModcom(A)$, \ie, $\dgMod(A) = Z^0(\dgModcom(A))$.
\edf

We now describe the projective and injective model structures on $\dgMod(A)$. For any $M_j \in \dgMod(A_j)$, define $\lambda_j(M_j) \in \dgMod(A)$ by
\begin{equation*}  
		(\lambda_j(M_j))_i \, := \, \bigoplus_{\varphi \in \Hom_I(j,i)} M_j \otimes_{A_j} A_i
\end{equation*}
with the obvious transition maps for $i \ra i'$ in $I$.
Dually, define $\rho_j(M_j) \in \dgMod(A)$ by
\begin{equation*}
		(\rho_j(M_j))_i \, := \, \prod_{\varphi \in \Hom_I(i,j)} M_j
\end{equation*}
with the obvious induced $A_i$-module structure, and with the obvious transition maps for $i \ra i'$ in $I$.

If $U_i : \dgMod(A) \ra \dgMod(A_i)$ is the forgetful functor $M \mapsto M_i$, then there is a three-way adjunction $\lambda_i \dashv U_i \dashv \rho_i$. These can be combined to a three-way adjunctions $\lambda \dashv U \dashv \rho$ between the categories
\begin{equation*}
	\begin{tikzcd}
	 \prod_{\alpha \in \Ob(I)} \, \dgMod(A_i)  \ar[r, shift left = 2] \ar[r, shift right = 2]
		& \dgMod(A)  \ar[l]
	\end{tikzcd}
\end{equation*}
where $U(M) = (M_i)_{i \in \Ob(I)}$, with left adjoint $\lambda$ given by $\lambda((M_i)_{i \in \Ob(I)}) = \bigoplus_{i \in \Ob(I)} \lambda_i(M_i)$, and right adjoint $\rho$ given by  $\rho((M_i)_{i \in \Ob(I)}) = \prod_{i \in \Ob(I)} \rho_i(M_i)$.

Applying the main results of \cite{HKRS17, GKR20}, we can use the pair $\lambda \dashv U$ and $U \dashv \rho$ to lift the product model structure on $\prod_{\alpha \in \Ob(I)} \, \dgMod(A_i)$ to $\dgMod(A)$ in two ways. Namely, in the \emph{projective model structure}, a map $f : M \ra N$ in $\dgMod(A)$ is a weak equivalence or fibration if and only if each $f_i : M_i \ra N_i$ is a weak equivalence or fibration in $\dgMod(A_i)$. 
In the \emph{injective model structure}, a map $f : M \ra N$ in $\dgMod(A)$ is a weak equivalence or cofibration if and only if each $f_i : M_i \ra N_i$ is a weak equivalence or cofibration in $\dgMod(A_i)$.

The DG enrichment $\dgModcom(A)$ of $\dgMod(A)$ can be described alternatively by either the copowering or the powering of $\Ch(k)$ on $\dgMod(A)$, both of which are defined pointwise, in the standard way. From this, it is easy to see that both the projective and injective structures are compatible with these powering/copowerings (for example, check copowering for injective model structure and powering for projective model structure), so that both model structures make $\dgMod(A)$ into a $\Ch(k)$-model structure in the sense of \cite[Definition 4.2.18]{Hov99}. This allows one to define the derived Hom complex ${\bm R}\Homcom_A(P,M)$ for $P,M \in \dgMod(A)$.

We now describe a standard cobar construction that computes this derived Hom complex. For any $P,M \in \dgMod(A)$, we consider the cosimplicial cochain complex $\Homcom_A^{\hc}(P,M)^{\bullet} \in \Ch(k)^{\Delta}$ by 
\begin{equation}  \label{Hom_hc_def}
\Homcom_A^{\hc}(P,M)^n \, := \, \prod_{F:[n] \ra I} \, \Homcom_{A_{F(0)}}(P_{F(0)},M_{F(n)})
\end{equation}
with the obvious cosimplicial structure maps between them.
The complex $\Homcom_A^{\hc}(P,M)$ is defined as the homotopy limit
\begin{equation*}
\Homcom_A^{\hc}(P,M) \, := \, \Tot^{\Pi}_{\Delta} \, \, \Homcom_A^{\hc}(P,M)^{\bullet}
\end{equation*}

Notice that the (strict) limit $\Homcom_A(P,M) := \ker(d^0 - d^1)$ of $\Hom_A^{\hc}(P,M)^{\bullet}$ is precisely the Hom complex from $P$ to $M$ in the DG category $\dgModcom(A)$. We think of the homotopy limit $\Homcom_A^{\hc}(P,M)$ as ``homotopy coherent'' maps from $P$ to $M$, hence the superscript ``hc''.

\bpp  \label{Homhc_RHom}
If each $P_i$ is cofibrant in $\dgMod(A_i)$, then the complex
$\Homcom_A^{\hc}(P,M)$ is quasi-isomorphic to the derived Hom complex ${\bm R}\Homcom_A(P,M)$.
\epp

\bpf
For any $A$-modules $P,M$, let $\cR^{\bullet} M \in \dgMod(A)^{\Delta}$ and $\cQ_{\bullet} P \in \dgMod(A)^{\Delta^{\op}}$ be defined by
\begin{equation*}
\begin{split}
(\cR^{n} M)_{i} \, &:= \, \prod_{F: [n] \ra I} \, \rho_{F(0)}(M_{F(n)}|_{F(0)}) \, = \, \prod_{F : [n] \ra I} \,\, \prod_{\varphi \in \Hom_I(i,F(0))} \, M_{F(n)} \\
(\cQ_{n} P)_{i} \, &:= \, \bigoplus_{F: [n] \ra I} \, \lambda_{F(n)}( P_{F(0)} \otimes_{A_{F(0)}} A_{F(n)} ) \, = \, \bigoplus_{F : [n] \ra I} \,\, \bigoplus_{\varphi \in \, \Hom_I(F(n),i)} P_{F(0)} \otimes_{A_{F(0)}} A_i
\end{split}
\end{equation*}
with the obvious $A$-module structure and (co)simplicial maps.
Then we have canonical isomorphisms of cosimplicial cochain complexes
\begin{equation*}
\Homcom_A(P,\cR^{\bullet}(M)) \, \cong \, \Homcom_A^{\hc}(P,M)^{\bullet}  \, \cong \, \Homcom_A(\cQ_{\bullet}(P),M) 
\end{equation*}

There are natural maps $M \ra \cR^0 M$ making $\cR^{\bullet} M$ an augmented cosimplicial $A$-module. Moreover, for each fixed $i \in I$, the underlying augmented system $M_i \ra (\cR^{\bullet} M)_i$ of $A_i$-modules has an extra degeneracy by virtue of the fact that $(\varphi,F)$ together forms a string of length $n+1$ in $I$. Thus, if we let $\cR(M) := \Tot^{\Pi}_{\Delta} \, \cR^{\bullet} M$ then there is a quasi-isomorphism $r_M : M \ra \cR (M)$ in $\dgMod(A)$. 
Dually, let $\cQ(P) := \Tot^{\oplus}_{\Delta^{\op}} \, \cQ_{\bullet} P$, then there is a natural quasi-isomorphism $q_P : \cQ(P) \ra P$ in $\dgMod(A)$. 


From the definition, it is easy to see that $\cR(M)$ is fibrant in the injective model structure, and $\cQ(P)$ is cofibrant in the projective model structure if each $P_i$ is cofibrant. Hence, $\Homcom_A^{\hc}(P,M)^{\bullet}$ computes the derived Hom complex in either model structure.
%
%
\epf


\bdf
Given a small category $I$, and a functor $A : I \ra \DGA_k$.
A DG module $M$ over $A$ is said to be \emph{strictly Cartesian} if the maps $M_i \otimes_{A_i} A_j \ra M_j$ is an isomorphism for each $i \ra j$ in $\Delta$. It is said to be \emph{homotopy Cartesian} if the maps $M_i \otimes_{A_i}^{{\bm L}} A_j \ra M_j$ is a quasi-isomorphism for each $i \ra j$ in $\Delta$.
\edf

\brm  \label{check_at_di}
We will mostly be interested in the case $I = \Delta$. Given $A \in \DGA_k^{\Delta}$, since $s^i \circ d^i = \id$ in $\Delta$, it suffices to check the Cartesian and homotopy Cartesian property of $M \in \dgMod(A)$ along $d^i : [n] \ra [n+1]$. In our main examples, the maps $d^i : A_n \ra A_{n+1}$ are flat, hence if $M$ is strictly Cartesian, it is also homotopy Cartesian.
\erm

Given $P,M \in \dgMod(A)$, suppose that $P$ is strictly Cartesian. Then one can define
\begin{equation}  \label{Hom_sim_def}
\Homcom_A^{\sim}(P,M) \, : \, I \ra \Ch(k) \, , \qquad \text{given by } \quad \Homcom_A^{\sim}(P,M)_i := \Hom_{A_i}(P_i,M_i)
\end{equation}
with transition map for $\varphi \in \Hom_I(i,j)$ given by
\begin{equation}   \label{Hom_cospan}
\Homcom_{A_i}(P_i,M_i) \raq \Homcom_{A_i}(P_i,M_j) \xlaq{\cong} 
\Homcom_{A_j}(P_j,M_j)
\end{equation}

Notice also that, since $P$ is strictly Cartesian, $\Homcom_A^{\hc}(P,M)^{\bullet}$ can be alternatively written as
\begin{equation}  \label{Hom_hc_cart}
\Homcom_A^{\hc}(P,M)^n \, \cong \, \prod_{F:[n] \ra I} \, \Homcom_{A_{F(n)}}(P_{F(n)},M_{F(n)})
\end{equation}

Denote by $\El(K) := \Delta \downarrow K$ the category of elements of a simplicial set $K$. Apply this to the nerve $N(I)$ of $I$, then we have functors
\begin{equation*}
\Delta \xlaq{p} \El(N(I)) \xraq{q} I
\end{equation*}
where $q$ sends $F : [n] \ra I$ to $F(n)$. Then the observation \eqref{Hom_hc_cart} can be rephrased by saying that 
\begin{equation*}
\Homcom_A^{\hc}(P,M) \, \cong \, \Ran_p(\, q^* \Homcom_A^{\sim}(P,M) \,)
\end{equation*}
As in, \eg, \cite[Section 4.4]{Rie}, since ${\bm R}\Ran_p \simeq \Ran_p$ and since the functor $q$ is homotopy initial%
\footnote{For any $\alpha \in \Ob(I)$, consider the slice category $q/\alpha$, whose objects consist of $F : [n] \ra I$ together with a map $\varphi : F(n) \ra \alpha$ in $I$, and whose morphisms consist of maps $[n] \ra [m]$ in $\Delta$ that respect both structures. Consider the endofunctor $S : q/\alpha \ra q/\alpha$ that sends $(F,\varphi)$ to $(\widetilde{F},\id)$, where $\widetilde{F} : [n+1]\ra I$ is obtained by combining $F$ and $\varphi$. Consider the object $* \in \Ob(q/\alpha)$ given by $F : [0] \ra I$, $F(0) = \alpha$, together with $\id : \alpha \ra \alpha$. Write also by $*$ the endofunctor of $q/\alpha$ that sends every object to $*$ and every morphism to the identity, then there are natural transformations $\id \Rightarrow S \Leftarrow *$ given by the maps $[n] \ra [n+1] \la [0]$ into the first $n+1$ elements and the last element. Since natural transformation become homotopy after taking nerves, we see that $N(q/\alpha)$ is contractible}, this shows that
\begin{equation}  \label{Homhc_holim_Homsim}
\parbox{40em}{ $\Homcom_A^{\hc}(P,M)$ is quasi-isomorphic to $\holim_{I} \, \Homcom_A^{\sim}(P,M)$.}
\end{equation}

More generally, if $P$ is homotopy Cartesian, and we take ${\bm R}\Hom$ everywhere, then the left pointing map in \eqref{Hom_cospan} is only a quasi-isomorphism. In this case, the formula ${\bm R}\Homcom_A^{\sim}(P,M)_i := {\bm R}\Homcom_{A_i}(P_i,M_i)$ gives an $\infty$-functor $I \ra {\it Ch}(k)$, and we have
\begin{equation}  \label{RHom_inflim_Homsim}
\parbox{40em}{ ${\bm R}\Homcom_A(P,M)$ is quasi-isomorphic to ${\rm lim}^{\infty}_{I} \, {\bm R}\Homcom_A^{\sim}(P,M)$.}
\end{equation}

\vspace{0.2cm}

Now assume that each $A_i$ is a non-positively graded CDGA. Let $\frX = \hocolim_{I^{\op}} \, \Spec A_i$, taken in the $\infty$-category (or model category) of derived stacks. Then the $\infty$-category of quasi-coherent sheaves on $\frX$ may be defined as the homotopy limit $\holim_{i \in I} \, {\it dgMod}(A_i)$, where we write ${\it dgMod(A_i)}$ for the $\infty$-category obtained by localizing the ordinary category $\dgMod(A_i)$ at the quasi-isomorphisms. In \cite{Ber12} (see also \cite{Pri13}), it is shown that this can be described as the $\infty$-category ${\it dgMod}(A)_{\hca}$ obtained by localizing the category $\dgMod(A)_{\hca}$ of homotopy Cartesian $A$-modules at the (pointwise) quasi-isomorphisms.

\vspace{0.2cm}

We now investigate our main case of interest: that of a global quotient. Let $G$ be a linear algebraic group acting on a non-positively graded CDGA $B$, and let $X = [Y/G]$ for $Y = \Spec \, B$ (see \eqref{Y_mod_G_def} for our convention). Let $A^{\bullet} = \cO(X_{\bullet})$. In this case, a strictly Cartesian $A^{\bullet}$-module $M^{\bullet}$ is simply a $G$-equivariant $B$-module $M^0$. Namely, the transition map 
$(d_0)^* M^0 \xra{\cong} M^1 \xla{\cong} (d_1)^* M^0$ gives a $G$-equivariant structure on $M^0$.
Conversely, given a $G$-equivariant $B$-module $M$ with $\cO(G)$-comodule structure map $\Delta : M \ra M \otimes \cO(G)$ (linear over $\Delta : B \ra B \otimes \cO(G)$), one can take
\begin{equation}  \label{assoc_cart}
C^n(G,M) \, := \, M \otimes \cO(G)^{\otimes n}
\end{equation}
with cosimplicial maps defined by a formula dual to \eqref{Y_mod_G_def}, making it into a strictly Cartesian $A^{\bullet}$-module. 
These two process are inverse to each other, giving an equivalence of categories $\dgMod_G(B) \simeq \dgMod(A)_{\ca}$, where $\dgMod(A)_{\ca} \subset \dgMod(A)$ is the full subcategory consisting of strictly Cartesian $A$-modules.
Given $G$-equivariant $B$-modules $M$ and $N$, we may form the $A^{\bullet}$-module $\Homcom_{A^{\bullet}}^{\sim}(C^{\bullet}(G,M), C^{\bullet}(G,N))$ as in \eqref{Hom_sim_def}, which will simply be denoted as $C^{\bullet}(G,\Homcom_B(M, N))$. In other words, we have
\begin{equation}  \label{Hom_sim}
C^n(G,\Homcom_B(M, N)) \, = \, \Homcom_B(M, N\otimes \cO(G)^{\otimes n})
\end{equation}
where the $B$-module structure on $N\otimes \cO(G)^{\otimes n}$ is induced by the map $\id_B \otimes 1 : B \ra B \otimes \cO(G)^{\otimes n}$ (see the proof of Lemma \ref{CG_Hom_RG} below for the details of the identification \eqref{Hom_sim}). 
From the above arguments (see \eqref{Homhc_holim_Homsim}), we see that if $M$ is cofibrant in $\dgMod(B)$, then the derived Hom complex between $C^{\bullet}(G,M)$ and $C^{\bullet}(G,N)$ is computed as the direct product total complex of $C^{\bullet}(G,\Homcom_B(M, N))$.

If $M$ is finitely generated and projective as a graded module, then we may move the operation $-\otimes \cO(G)^{\otimes n}$ outside of $\Homcom_B(M,-)$, so that $C^{\bullet}(G,\Homcom_B(M, N))$ is strictly Cartesian. This defines the structure of a $G$-equivariant module on $\Homcom_B(M, N)$, so that $C^{\bullet}(G,\Homcom_B(M, N))$ is the associated $A$-module. In general, however, \eqref{Hom_sim} is a slight abuse of notation, and is not intended to mean that $\Homcom_B(M, N)$ is an equivariant module in the strict sense (although we would like to think of the existence of the $\cO(X)^{\bullet}$-module structure on \eqref{Hom_sim} as an equivariant structure of some sort on $\Homcom_B(M, N)$).

There is an alternative description of $C^{\bullet}(G,\Homcom_B(M, N))$. 
Recall that the forgetful functor $U : \dgMod_{G}(B) \ra \dgMod(B)$ has a right adjoint $R$ given by $R(N) = N \otimes \cO(G)$, with $B$-module structure induced by the map $\Delta : B \ra B \otimes \cO(G)$. The Barr-Beck construction then associates to each $M \in \dgMod_{G}(B)$ an augmented cosimplicial object $M \ra \cR_G^{\bullet}(M)$ in $\dgMod_{G}(B)$, where $\cR_G^{n}(M) = M \otimes \cO(G)^{\otimes n+1}$.
The cosimplicial complex $C^{\bullet}(G,\Homcom_B(M, N))$ may be alternatively described as follows:
\blm  \label{CG_Hom_RG}
There is an isomorphism of cosimplicial complex 
$C^{\bullet}(G,\Homcom_B(M, N)) \cong \Homcom_B^G(M, \cR_G^{\bullet}(M) )$,
where we denote by $\Homcom_B^G(-,-)$ the Hom complex in the DG category $\dgModcom_G(B)$.
\elm

\bpf
By definition, $C^{\bullet}(G,\Homcom_B(M, N)) := \Homcom_{A^{\bullet}}^{\sim}(C^{\bullet}(G,M), C^{\bullet}(G,N))$ is given in degree $n$ by 
\begin{equation*}
	C^{n}(G,\Homcom_B(M, N)) \, = \, \Homcom_{B \otimes \cO(G)^{\otimes n}}( M \otimes \cO(G)^{\otimes n}, N \otimes \cO(G)^{\otimes n} )
\end{equation*}
There are two useful conventions to simplify this expression. Namely, each map $\varphi : [0] \ra [n]$ in $\Delta$ gives rise to a map $\varphi_* : M \ra M \otimes \cO(G)^{\otimes n}$ linear and Cartesian over the map $\varphi_* : B \ra B \otimes \cO(G)^{\otimes n}$ of CDGAs. We may fix a choice of such $\varphi$ for each $n$, and use the induction-restriction adjunction to simplify the expression.

If we choose $\varphi_0 : [0] \ra [n]$ by $\varphi_0(0) = 0$ for each $n$, then we obtain $C^{n}(G,\Homcom_B(M, N)) \cong \Homcom_B(M, N \otimes \cO(G)^{\otimes n})$ with the $B$-module structure on $N \otimes \cO(G)^{\otimes n}$ induced via $(\varphi_0)_* = \id_B \otimes 1 : B \ra B \otimes \cO(G)^{\otimes n}$. This description is precisely the one used in \eqref{Hom_sim}. One can verify that the induced cosimplicial structure maps are precisely those induced from a $G$-equivariant $B$-module structure on $\Homcom_B(M,N)$ if $M$ is finitely generated and projective as a graded $B$-module.

If we choose $\varphi_n : [0] \ra [n]$ by $\varphi_0(0) = n$ for each $n$, then we obtain $C^{n}(G,\Homcom_B(M, N)) \cong \Homcom_B(M, N \otimes \cO(G)^{\otimes n})$ with the $B$-module structure on $N \otimes \cO(G)^{\otimes n}$ induced via $(\varphi_n)_* = \Delta^{(n)} : B \ra B \otimes \cO(G)^{\otimes n}$. This is therefore identified with 
\begin{equation*}
	\Homcom_B^G(M, \cR_G^{n}(N) ) \, = \, \Homcom_B^G(M, R(R^n(N)) ) \, = \, \Homcom_B( M , R^n(N) )
\end{equation*}
One can verify that the cosimplicial structure maps are also the same.
\epf

Consider again the adjunction $U\dashv R$ between $\dgMod(B)$ and $\dgMod_G(B)$. For each $M \in \dgMod_G(B)$, the adjunction unit is given by 
the $G$-equivariant structure map $\Delta_M : M \ra M \otimes \cO(G)$. One can ask whether this map splits: 

\bdf
A $G$-equivariant DG module $M$ over $B$ is said to have \emph{split diagonal} if the map $\Delta_M : M \ra M \otimes \cO(G)$ has a left inverse in $\dgMod_G(B)$. 
\edf

If $M \in \dgMod_G(B)$ has split diagonal, then the splitting map induces in the obvious way an extra degeneracy $s^{-1}$ for the augmented cosimplicial object $M \ra \cR_G^{\bullet}(M)$ in $\dgMod_{G}(B)$. Hence, Lemma \ref{CG_Hom_RG} implies the following
\bcor  \label{split_diag_equiv}
If $N \in \dgMod_{G}(B)$ has split diagonal, then the canonical map
\begin{equation}  \label{Hom_G_to_CG_Hom}
	\Homcom_B^G(M,N) \raq \Tot^{\Pi}_{\Delta} \, C^n(G,\Homcom_B(M, N))
\end{equation}
is a quasi-isomorphism.
\ecor

Now we consider (linearly) reductive groups. We will show that if $G$ is (linearly) reductive, then every object in $\dgMod_{G}(B)$ has split diagonal (see Corollary \ref{reductive_split_diag} below). First, we recall the existence of Reynold operators on Hom complexes%
\footnote{We learned about this (for $B = k$) from \cite{CHK05}, but it should be quite standard.}:
\bpp  \label{Reynold_prop}
For each $M,N \in \dgMod_{G}(B)$, there is a canonical retract $E : \Homcom_{B}(M,N) \ra \Homcom_{B}^G(M,N)$. Moreover, for any $K,L \in \dgMod_{G}(B)$, the following diagram commutes:
\begin{equation}  \label{Reynold_comm}
	\begin{tikzcd}
		\Homcom_B^G(N,L) \otimes \Homcom_{B}(M,N) \otimes \Homcom_B^G(K,M) \ar[r, "\circ"] \ar[d, "\id \otimes E \otimes \id"'] 
		& \Homcom_{B}(K,L) \ar[d, "E"]\\
		\Homcom_B^G(N,L) \otimes \Homcom_{B}^G(M,N) \otimes \Homcom_B^G(K,M) \ar[r, "\circ"] 
		& \Homcom_{B}^G(K,L)
	\end{tikzcd}
\end{equation}
\epp

\bpf
When $B = k$ and when $K,L,M,N$ are finite dimensional vector spaces concentrated in degree $0$, then $E$ is the standard Reynold operator (\ie, the unique $G$-equivariant retract $\Hom_k(M,N) \ra \Hom_k^G(M,N)$), and \eqref{Reynold_comm} is a basic property. When $M,N$ are not necessarily finite dimensional (but still concentrated in degree $0$, with $B = k$), we can characterize the map $E$ as follows:
\begin{equation}  \label{Reynold_1}
	\parbox{40em}{Suppose $f : M \ra N$ satisfies $f(V) \subset W$ for finite dimensional $G$-stable subspaces $V \subset M$ and $W \subset N$, denote by $f_V^W : V \ra W$ the restriction, then $E(f)|_V = E(f_V^W)$.}
\end{equation}
If $M,N$ are finite dimensional, then \eqref{Reynold_1} is a property of $E$, which follows from the finite dimensional case of \eqref{Reynold_comm}.
Since every element in a rational $G$-module is contained in a finite dimensional $G$-stable subspace, the property \eqref{Reynold_1} uniquely specify $E(f)$. Namely, for each $x \in M$, choose such $V,W$ such that $x \in V$, so that \eqref{Reynold_1} determines $E(f)(x)$. The property \eqref{Reynold_1} in the finite dimensional case can in turn be used to show that $E(f)(x)$ is independent of the choice of $V,W$ so that $E$ is well-defined. It clearly satisfies \eqref{Reynold_comm}.

Now if $B = k$, and $M,N$ are any cochain complex of rational $G$-modules, then one can apply $E$ degreewise to $\Homcom^i_k(M,N) = \prod_{j \in \bZ}\Hom_k(M^j,N^{i+j})$. It suffices to verify that the obtained retract $E : \Homcom_{k}(M,N) \ra \Homcom_{k}^G(M,N)$ is a map of cochain complexes. This in turn follows from \eqref{Reynold_comm} in the above case because the differentials $d : M^i \ra M^{i+1}$ (and similarly for $N$) are $G$-equivariant maps.

Finally, consider the general case of $B$ being any CDGA. Since $\Homcom_{B}(M,N) \subset \Homcom_{k}(M,N)$ is a subcomplex, it suffices to show that $E$ preserves this subcomplex. Notice that the fact that $\Delta_M : M \ra M \otimes \cO(G)$ is linear over $\Delta_B : B \ra B \otimes \cO(G)$ can be rewritten as the fact that $\mu_M : B \otimes M \ra M$ is $G$-linear. Thus, given any $f \in \Homcom_B(M,N)$, we may apply $E$ to $f \circ \mu_M = \mu_N \circ (f \otimes \id_B)$, so that by applying \eqref{Reynold_comm} for the case of $B = k$, we have
\begin{equation*}
	E(f) \circ \mu_M = E(f \circ \mu_M) = E(\mu_N \circ (f \otimes \id_B)) = \mu_N \circ E(f \otimes \id_B) = \mu_N \circ E(f) \otimes \id_B
\end{equation*}
where the last equality follows from the following fact:
\begin{equation*}
	\parbox{40em}{For any $f \in \Homcom_k(M,N)$ and $g \in \Homcom_k^G(K,L)$, denote their tensor product by $f \otimes g \in \Homcom_k(M \otimes K,N \otimes L)$, then we have $E(f \otimes g) = E(f) \otimes g$.}
\end{equation*}
which is obvious in the finite dimensional case, which then implies the infinite dimensional case, and then the case for cochain complexes.
\epf

\bcor  \label{reductive_split_diag}
If $G$ is (linearly) reductive, then every $N \in \dgMod_{G}(B)$ has split diagonal, hence the map \eqref{Hom_G_to_CG_Hom} is a quasi-isomorphism.
\ecor

\bpf
Notice that the map $f = \id_N \otimes \epsilon: N \otimes \cO(G) \ra N$ is $B$-linear (where $N \otimes \cO(G)$ has the $B$-module structure induced by $\Delta : B \ra B \otimes \cO(G)$). Hence $E(f)$ is a map in $\dgMod_G(B)$. Moreover, since $f \circ \Delta_N = \id$, the property \eqref{Reynold_comm} then implies that $E(f) \circ \Delta_N = \id$, so that $N$ has split diagonal.
\epf

\brm
For $G$ (linearly) reductive, $B = k$, $N$ bounded below, and $M$ bounded and finite dimensional, the fact that \eqref{Hom_G_to_CG_Hom} is a quasi-isomorphism can be proved by the standard methods in classical homological algebra.
\erm

\brm
The functor $C^{\bullet}(G,-) : \dgMod_G(B) \ra \dgMod(A^{\bullet})$ preserves quasi-isomorphisms, and hence induces an $\infty$-functor $dgMod_G(B) \ra dgMod(A^{\bullet})$ between the $\infty$-categories obtained by localizing at quasi-isomorphisms. This functor lands in the full subcategory $dgMod(A^{\bullet})_{\hca} \subset dgMod(A^{\bullet})$ consisting of homotopy Carteisan modules. We expect that the $\infty$-functor $dgMod_G(B) \ra dgMod(A^{\bullet})_{\hca}$ is an equivalence for any linear algebraic group $G$, but have not been able to prove it. Since we will work entirely within $dgMod(A^{\bullet})_{\hca}$, we will not need this relation in this paper.
\erm

\section{Small cosimplicial algebras}  \label{app_small_alg}

In the first part of this Appendix, we identify a certain subcategory of cosimplicial commutative algebras that is equivalent to the category of CDGAs generated in degrees $0$ and $1$. Most of our arguments in this part are parallel to those in \cite{Rap88}, except that we do not impose the reducedness assumption. We supply the details whenever modifications are needed for this generality. 
After that, we discuss some applications to equivariant modules.

Given a cosimplicial $k$-module $A$, denote by $C(A)$ its total complex $C^n(A) = A^n$, with differential $d = \sum_{i=0}^{n+1} (-1)^i {d^i}$. Denote by $N(A)$ the normalized subcomplex, defined by $N^n(A) = \cap_{i = 0}^{n-1} \ker(s^i : A^n \ra A^{n-1})$, which is a subcomplex of $(C(A),d)$. Recall that $D^n(A) := \sum_{i = 1}^{n} {\rm Im}(d^i : A^{n+1} \ra A^n)$ is also a subcomplex, which forms a direct sum $C(A) = N(A) \oplus D(A)$. The subcomplex $D(A)$ has zero cohomology, so that the inclusion $N(A) \rinto C(A)$ is a quasi-isomorphism. 
Denote the cosimplicial Dold-Kan correspondence by
\begin{equation*}
\begin{tikzcd}
N \, : \, \Mod_k^{\Delta} \ar[r, shift left]
& \dgMod_k^{\geq 0} \, : \, \Gamma \ar[l, shift left]
\end{tikzcd}
\end{equation*}
For any $A,B \in \Mod_k^{\Delta}$, there are the Alexander-Whitney maps and the Eilenberg-Zilber shuffle maps
\begin{equation*}
\AW : C(A) \otimes C(B) \ra C(A \otimes B) \qquad \text{and} \qquad 
\nabla : C(A \otimes B) \ra  C(A) \otimes C(B)
\end{equation*}
both of which preserve the normalized subcomplexes. \ie, they restrict to  
\begin{equation*}
\AW : N(A) \otimes N(B) \ra N(A \otimes B) \qquad \text{and} \qquad 
\nabla : N(A \otimes B) \ra  N(A) \otimes N(B)
\end{equation*}
Moreover, $\nabla$ is a left inverse to $AW$ when restricted to the normalized subcomplexes. \ie, we have $\nabla \circ \AW = \id$ on $N(A) \otimes N(B)$. 

The Alexander-Whitney maps are associative, while the Eilenberg-Zilber maps are associative and commutative. In particular, if $A \in \Alg_k^{\Delta}$, then $C(A)$ is a DGA, with a cup product $\cup$ induced by the Alexander-Whitney map. Moreover, $N(A)$ is a sub-DGA. Explicitly, the cup product is given by
\begin{equation*}
x \cup y \, := \, (d^{p+q}\ldots d^{p+1}(x) )\cdot (d^0 \ldots d^0(y))
\qquad \quad \text{for } x \in A^p, \, \, y \in A^q
\end{equation*}

If confusion is possible, we will call the original product $A^n \otimes A^n \ra A^n$ the \emph{degreewise product}, and write a dot for the product, to distinguish it from the cup product $\cup : A^p \otimes A^q \ra A^{p+q}$. 
It is clear that we have
\begin{equation}  \label{cup_d}
\begin{split}
d^i(x) \cup y &= d^i(x \cup y) \qquad \text{for } 0 \leq i \leq p \\
x \cup d^j(y) &= d^{j+p}(x \cup y) \qquad \text{for } 1 \leq j \leq q+1
\end{split}
\end{equation}
from which it follows that $DA \subset CA$ is a left ideal%
\footnote{There seems to be a small mistake in \cite{Rap88}. Namely, the subcomplex $DA$ is defined in \cite{Rap88} using a different convention, which seems to be then a \emph{right} ideal.} under the cup product.

Similarly, the Eilenberg-Zilber maps can be rephrased as a lax monoidal structure on $\Gamma$, so that if $B \in \DGAp_k$, then $\Gamma(B)$ has a shuffle product. If $B \in \CDGAp_k$, then the shuffle product is moreover commutative. This gives two functors
\begin{equation*}
N : \Alg_k^{\Delta} \ra \DGAp_k \, , \qquad \quad  \text{and} \qquad \quad  
\Gamma : \DGAp_k \ra \Alg_k^{\Delta}
\end{equation*}
We have the following
\blm  \label{Gamma_ff}
The composition $N \circ \Gamma$ is isomorphic to the identity functor on $\DGAp_k$. Moreover, the functor $\Gamma$ is always fully faithful.
\elm

\bpf
The first statement follows from the fact that $\nabla \circ \AW = \id$ on normalized subcomplexes. The second statement then follows by considering the commutative diagram
\begin{equation*}
\begin{tikzcd}
\Hom_{\DGAp_k} (A,B) \ar[r, "\Gamma"] \ar[rr, bend left = 10, "\cong"] \ar[d, hook]
&\Hom_{\Alg_k^{\Delta}}(\Gamma(A) , \Gamma(B)) \ar[r, "N"]  \ar[d, hook]
& \Hom_{\DGAp_k} (N\Gamma(A),N\Gamma(B))  \ar[d, hook] \\
\Hom_{\dgMod^{\geq 0}_k} (A,B) \ar[r, "\Gamma", "\cong"']
&\Hom_{\Mod_k^{\Delta}}(\Gamma(A) , \Gamma(B)) \ar[r, "N", "\cong"']
& \Hom_{\dgMod^{\geq 0}_k} (N\Gamma(A),N\Gamma(B)) 
\end{tikzcd}
\end{equation*}
\epf

Given a class of non-negatively graded DGAs, the functor $\Gamma$ then embeds it as a full subcategory of $\Alg^{\Delta}_k$. We may then try to characterize this subcategory by specifying some properties that these objects should satsify. We will focus primarily on CDGAs. Then $\Gamma(B)$ must be commutative. Moreover, it follows directly from the definition of the shuffle product that $\Gamma(B)$ must be $1$-thin in the sense of the following

\bdf
A cosimplicial commutative algebra $A \in \Alg_k^{\Delta}$ is said to be \emph{$1$-thin} if the ideal $N^1(A) = \ker(s^0) \subset A^1$ is square zero (with respect to the degreewise product).
\edf

In fact, if $B$ is a non-negatively graded CDGA, then each $\Gamma^n(B)$ is an infinitesimal thickening of $B^0$, and there are many more ``thinness'' conditions that $\Gamma(B)$ should specify. However, for our purpose, this condition suffices.

\bdf
A simplicial algebra $A \in \Alg_k^{\Delta}$ is said to be $[0,1]$-generated if the DGA $N(A)$ is generated in degree $0$ and $1$ (by cup product and $k$-linear sums).

A simplicial algebra $A$ is said to be \emph{small} if it is commutative, $1$-thin, and $[0,1]$-generated.

A non-negatively graded CDGA $B$ is said to be \emph{small} if it is generated in degree $0$ and $1$.
\edf

\blm  \label{NA_01gen_CA}
If $A \in \Alg_k^{\Delta}$ is $[0,1]$-generated, then $C(A)$ is generated in degree $0$ and $1$ (by cup product and $k$-linear sums).
\elm

\bpf
By the decomposition $C^n(A) = N^n(A) \oplus D^n(A)$, it suffices to show that $D^n(A)$ is generated by cup product by $C^1(A)$, for each $n \geq 2$.
We first prove it for $n=2$. For each $x \in C^1(A)$, notice that $d^0(x) = 1_{A_1} \cup x$ and $d^2(x) =  x \cup 1_{A_1}$, hence both are $\cup$-generated by $C^1(A)$. If $x \in N^1(A)$, then $d^0(x) - d^1(x) + d^2(x) \in N^2(A)$, which is $\cup$-generated by $N^1(A)$ by assumption. Thus, $d^1(x)$ is also $\cup$-generated by $C^1(A)$. Finally, if $x \in D^1(A)$, say $x = d^1(y)$, then we have $d^1(x) = d^2(x)$, hence also $\cup$-generated by $C^1(A)$.
To complete the proof, one can prove the case for $n \geq 3$ by induction, using \eqref{cup_d}.
%
\epf

\brm
In \cite{Rap88}, the converse is also shown to be true under a further reducedness assumption (\ie, $A_0 = k$). We haven't been able to prove the converse in general.
\erm

\blm
If $A \in \CAlg^{\Delta}_k$ is small, then $N(A)$ is commutative, hence small.
\elm

\bpf
Since $N(A)$ is generated in degrees $0$ and $1$, it suffices to show that these elements (graded) commute with each other. The fact that $N^0(A)$ is commutative is obvious. Given $a \in N^0(A)$ and $x \in N^1(A)$, we have $a \cup x = d^1(a) \cdot x$ and $x \cup a = x \cdot d^0(a)$. These are the same since $d^0(a) - d^1(a) \in N^1(A)$, which multiply with $x$ to get $0$. Thus, it suffices to show that $x \cup x = 0$, for which the proof of \cite[Lemma 2.8a]{Rap88} carries over without change.
\epf

Thus, we see that $\Gamma$ and $N$ give functors 
\begin{equation}  \label{N_Gamma_small}
\begin{tikzcd}
N \, : \, (\CAlg^{\Delta}_k)_{\rm small} \ar[r, shift left] &  (\CDGAp_k)_{\rm small} \, : \, \Gamma  \ar[l, shift left]
\end{tikzcd}
\end{equation}

\bpp
The functors \eqref{N_Gamma_small} are inverse equivalences of each other.
\epp

\bpf
It suffices to show that $A \cong \Gamma (N(A))$ for all $A \in (\CAlg^{\Delta}_k)_{\rm small}$. Clearly, they have the same simplicial module structure. Moreover, we have $N(A) \cong N(\Gamma(N(A)))$ in $\CDGAp_k$. Therefore, we may identify $A$ and $\Gamma (N(A))$ as a simplicial module, and think of having two commutative product structure $\cdot$ and $\cdot'$ on $A$ (both are small) inducing cup products $\cup$ and $\cup'$ on $C(A)$ respectively, such that $\cup = \cup'$ on $N(A) \subset C(A)$. 

First, we show that $\cup = \cup'$ on $C(A)$. Since both product structures are small, by Lemma \ref{NA_01gen_CA}, it suffices to show that 
\begin{equation}  \label{cup_a_equal}
a_1 \cup a_2 \cup \ldots \cup a_n \, = \, a_1 \cup' a_2 \cup' \ldots \cup' a_n
\end{equation}
for each $a_i$ that are in either $A^0$, $N^1(A)$ or $D^1(A)$. We induct on the number of times such elements come from $D^1(A)$. If none of the $a_i$'s are in $D^1(A)$, then \eqref{cup_a_equal} holds by assumption. If $a_i \in D^1(A)$ for some $i > 1$, say $a_i = d^1(b_i)$, then we can apply \eqref{cup_d} to both sides and use the induction hypothesis. If $a_1 = d^1(b)$, then notice that
\begin{equation*}
d^1(b) \cup y = d^0(b) \cup y + (d^0(b) - d^1(b)) \cup y = d^0(b \cup y) + (d^0(b) - d^1(b)) \cup y
\end{equation*}
where $y = a_2 \cup \ldots \cup a_n$ (the same holds for $\cup'$). Since $d^0(b) - d^1(b) \in N^1(A)$, we may use the induction hypothesis again.

Next, we note that $\cdot = \cdot'$ on $A^0$ and $A^1$. For $A^0$ it can be directly read off from the cup product. For $A^1$, since both products are commutative and square-zero on $N^1(A)$, it suffices to show that $d^1(a)  \cdot x = d^1(a) \cdot' x$ for $x \in N^1(A)$ and $a \in A^0$. But these are just $x \cup a$ and $x \cup' a$ respectively, and hence are equal.
To complete the proof (\ie, to show that $\cdot = \cdot'$), we note, as in \cite{Rap88}, that the commutativity of $A$ implies that
\begin{equation*}
(x_1 \cup \ldots \cup x_p) \cdot (y_1 \cup \ldots \cup y_p) = (x_1 \cdot y_1) \cup \ldots \cup (x_p \cdot y_p)
\end{equation*}
for $x_i, y_i \in A^1$ (the same holds for $\cdot'$ and $\cup'$).
\epf

\bdf
Given any $[0,1]$-generated simplicial commutative algebra $A$, its \emph{$1$-thin quotient} is the simplicial commutative algebra $A_{(1)}$ obtained by quotienting out the simplicial ideal generated by $N^1(A)^2 \subset A^1$. Clearly, $A_{(1)}$ is small.
\edf

Recall that we have a functor $\Gamma : \CDGA^{\geq 0}_k \ra \CAlg_k^{\Delta}$. It is easy to see that it has a left adjoint $\Gamma^* : \CAlg_k^{\Delta} \ra \CDGA^{\geq 0}_k$. Proposition \ref{glob_quot_1thin_quot} allows us to describe $\Gamma^*(A)$ when $A$ is $[0,1]$-generated:

\bpp  \label{Gamma_star_1thin}
Suppose that $A \in \CAlg_k^{\Delta}$ is $[0,1]$-generated, then we have $\Gamma^*(A) \cong N(A_{(1)})$, where $A_{(1)}$ is the $1$-thin quotient of $A$. Moreover, the adjunction counit $A \ra \Gamma \Gamma^*(A)$ is the quotient map $A \ra A_{(1)}$ to its $1$-thin quotient.
\epp

\bpf
For any $B \in \CDGA^{\geq 0}_k$, we have seen that $\Gamma(B)$ is $1$-thin, so that we have
$\Hom_{\CAlg_k^{\Delta}}(A,\Gamma(B)) \cong \Hom_{\CAlg_k^{\Delta}}(A_{(1)},\Gamma(B))$. Since $A_{(1)}$ is small, we have $A_{(1)} = \Gamma(N(A_{(1)}))$. Since $\Gamma$ is always fully faithful (see Lemma \ref{Gamma_ff}), we see that $\Hom_{\CAlg_k^{\Delta}}(A_{(1)},\Gamma(B)) \cong \Hom_{\CDGAp_k}(N(A_{(1)}),B)$, which shows that $\Gamma^*(A) \cong N(A_{(1)})$. The description of the adjunction counit as the quotient map is also clear by tracing through the isomorphisms.
\epf

We now consider the special case of a global quotient. Thus, let $G$ be a linear algebraic group acting on an affine scheme $Y = \Spec \, B$, and let $X = [Y/G]$ (see \eqref{Y_mod_G_def} for our convention).
In the cosimplicial commutative algebra given by $\cO(X)^n = B \otimes \cO(G)^{\otimes n}$. 
Its cup product is given by
\begin{equation*}  
(b' \otimes f_1 \otimes \ldots \otimes f_p) \cup (b'' \otimes f_{p+1} \otimes \ldots \otimes f_{n} ) \, = \, 
[(b' \otimes f_1 \otimes \ldots \otimes f_p) \cdot \Delta^{(p)}(b'')] \otimes f_{p+1} \otimes \ldots \otimes f_{n}
\end{equation*} 
where $\Delta^{(p)} : B \ra B \otimes \cO(G)^{\otimes p}$ is the repeated coaction map.
In particular, we have
$b \cup f_1 \cup \ldots \cup f_n = b \otimes  f_1 \otimes \ldots \otimes f_p$, where, on the left hand side, we think of $b$ to lie in $\cO(X_0)$ and $f_i$ to lie in $\cO(X_1)$ via $1 \otimes f_i$.

In a more classical language, we may write $C^n(G,B)$ as the set of algebraic functions $G^n \ra B$. Then for $F_1 \in C^p(G,B)$ and $F_2 \in C^q(G,B)$, the cup product is the function
\begin{equation}  \label{cup_OYG_2}
(F_1 \cup F_2)(g_1,\ldots,g_{p+q})
\, = \, 
F_1(g_1,\ldots,g_p) \cdot \rho_{g_1\ldots g_p}( F_2(g_{p+1},\ldots, g_{p+q}) )
\end{equation}
where $\rho$ is the left $G$-action on $B$ dual to the right $G$-action on $Y$.

Notice that the normalized subcomplex is given by $N^n(\cO(X)) = B \otimes \frm^{\otimes n}$, where $\frm \subset \cO(G)$ is the maximal ideal corresponding to the identity element $e \in G(k)$. Thus, $\cO(X)$ is $[0,1]$-generated.
We may then form its $1$-thin quotient, which is described by the following

\bpp  \label{glob_quot_1thin_quot}
The normalization of the $1$-thin quotient of $\cO([Y/G])$ is given by
$N( \cO([Y/G])_{(1)} ) \cong \CE(\frg^*,B)$ as a CDGA. 
If $M$ is a $G$-equivariant $B$-module, then $N( C^{\bullet}(G,M) \otimes_{\cO([Y/G])} \cO([Y/G])_{(1)} ) \cong \CE(\frg^*,M)$ as a DG module over $\CE(\frg^*,B)$.

Moreover, the map to these quotients are given by
\begin{equation}  \label{glob_quot_1thin_quot_map}
\begin{tikzcd}
N^n( \cO([Y/G]) ) \ar[r, twoheadrightarrow] \ar[d, equal]  &  N^n( \cO([Y/G])_{(1)} ) \ar[d, equal]  \\
B \otimes \frm^{\otimes n} \ar[r, twoheadrightarrow] 
& B \otimes \Lambda^n(\frg^*)
\end{tikzcd}
\end{equation}
\begin{equation}
\label{glob_quot_1thin_quot_map_module}
\begin{tikzcd}
N^n( C^{\bullet}(G,M) ) \ar[r, twoheadrightarrow] \ar[d, equal]  &  N^n( C^{\bullet}(G,M) \otimes_{\cO([Y/G])} \cO([Y/G])_{(1)} ) \ar[d, equal]  \\
M \otimes \frm^{\otimes n} \ar[r, twoheadrightarrow] 
& M \otimes \Lambda^n(\frg^*)
\end{tikzcd}
\end{equation}
where the horizontal map on bottom of each diagram is obtained by taking $\frm \ronto \frm/\frm^2 = \frg^*$ and then multiply.
\epp

\bpf
The statements for $B = k$ are precisely \cite[Lemma 3.1, 3.4]{Rap88}. In particular \cite[Lemma 3.4]{Rap88} establishes both of our statements at the level of cochain complexes, so we just need to check the compatibility of product and module structures, which is obvious.
\epf

A similar statement also holds for Hom complexes. Let $P,M \in \dgMod_G(B)$. Then recall from \eqref{Hom_sim} that there is an $\cO([Y/G])$-module $C^{\bullet}(G,\Homcom_B(P,M))$. Recall also from \eqref{Hom_g_comod} that there is a $\CE(\frg^*,B)$-module $\CE(\frg^*,\Homcom_B(P,M))$. 
These two are related by
\begin{equation}  \label{Hom_hG_to_CE}
\begin{tikzcd}
N^{n}(G,\Homcom_B(P,M)) \ar[r, twoheadrightarrow] \ar[d, equal]
& \CE^n(\frg^*, \Homcom_B(P,M)) \ar[d, equal]\\
\Homcom_B(P,M \otimes \frm^{\otimes n}) \ar[r, twoheadrightarrow]
& \Homcom_B(P,M \otimes \Lambda^n \frg^*)
\end{tikzcd}
\end{equation}
where the horizontal map on bottom is obtained by taking $\frm \ronto \frm/\frm^2 = \frg^*$ and then multiply.
This map is linear over the map \eqref{glob_quot_1thin_quot_map} of DG algebras.
Indeed, if $P$ is finitely generated and projective as a graded $B$-module, then $Q = \Homcom_B(P,M)$ has an honest equivariant $B$-module structure, so that this follows from \eqref{glob_quot_1thin_quot_map_module}. In general, the map on the first line of \eqref{Hom_hG_to_CE} is defined by the normalization of
\begin{equation}  \label{Hom_hG_to_CE_tensor}
\Homcom_{A^{\bullet}}^{\sim}( C^{
\bullet}(G,P), C^{
\bullet}(G,M) ) \raq 
\Homcom_{A^{\bullet}}^{\sim}( C^{
	\bullet}(G,P), C^{
	\bullet}(G,M) \otimes_{\cO([Y/G])} \cO([Y/G])_{(1)} )
\end{equation}
(see \eqref{Hom_sim_def} for notation) and the validity of the description \eqref{Hom_hG_to_CE} of it can be detected by taking maps $K \ra P$ from $K$ that is finitely generated and projective as a graded $B$-module.



\end{document}